\documentclass{amsart}
\usepackage{amssymb,amsmath,amscd, amsfonts,latexsym}
\usepackage[top = 4cm, bottom = 4cm, left = 4cm, right = 4cm]{geometry}
\usepackage{hyperref}
\usepackage{float}

\newtheorem{theorem}{Theorem}

\newtheorem{lemma}{Lemma}

\theoremstyle{remark}

\title[\title{The $x-$coordinates of Pell equations and sums of two Padovan numbers}]{On the $x-$coordinates of Pell equations which are sums of two Padovan numbers}

\author[Mahadi Ddamulira ]{Mahadi Ddamulira}



\address{Mahadi Ddamulira \newline
         \indent Institute of Analysis and Number Theory, Graz University of Technology \newline
        \indent Kopernikusgasse 24/II \newline
         \indent A-8010 Graz, Austria}

\email{mddamulira\char'100tugraz.at; mahadi\char'100aims.edu.gh}

\begin{document}


\begin{abstract}
Let $ \{P_{n}\}_{n\geq 0} $ be the sequence of Padovan numbers  defined by $ P_0=0 $, $ P_1 = P_2=1$ and $ P_{n+3}= P_{n+1} +P_n$ for all $ n\geq 0 $. In this paper, we find all positive square-free integers $ d $ such that the Pell equations $ x^2-dy^2 = \pm 1 $, $ X^2-dY^2=\pm 4 $ have at least two positive integer solutions $ (x,y) $ and $(x^{\prime}, y^{\prime})$,  $ (X,Y) $ and $(X^{\prime}, Y^{\prime})$, respectively, such that each of  $ x, ~x^{\prime}, ~X, ~X^{\prime} $ is  a sum of two Padovan numbers.
\end{abstract}

\maketitle

\section{Introduction}
\noindent
Let $ \{P_{n}\}_{n\geq 0} $ be the sequence of  Padovan numbers given by
$$P_0=0, ~P_1= 1, ~P_2=1, \text{  and  } P_{n+3}= P_{n+1}+P_n \text{   for all  } n\geq 0.$$
This is sequence $ A000931  $ on the Online Encyclopedia of Integer Sequences (OEIS).
The first few terms of this sequence are
$$\{P_{n}\}_{n\ge 0} = 0, 1, 1, 1, 2, 2, 3, 4, 5, 7, 9, 12, 16, 21, 28, 37, 49, 65, 86, 114, 151, \ldots. $$

\noindent
Let $d\geq 2$ be a positive integer which is not a square. It is well known that the Pell equations
\begin{eqnarray}
x^{2}-dy^{2}=\pm 1, \label{Pelleqn1}
\end{eqnarray}
and 
\begin{eqnarray}
X^{2}-dY^{2}=\pm 4, \label{Pelleqn2}
\end{eqnarray}
have infinitely many positive integer solutions $(x,y)$ and $ (X,Y) $, respectively. By putting $(x_1, y_1)$ and $ (X_1, Y_1) $ for the smallest positive solutions to \eqref{Pelleqn1} and \eqref{Pelleqn2}, respectively, all solutions are of the forms $ (x_k, y_k) $ and $ (X_k, Y_k) $ for some positive integer $k$, where
\begin{eqnarray*}
x_k+y_k\sqrt{d} = (x_1+y_1\sqrt{d})^k\qquad {\text{\rm for~all}} \quad k\ge 1,\label{Pellsoln1}
\end{eqnarray*}
and
\begin{eqnarray*}
\dfrac{X_k+Y_k\sqrt{d}}{2} = \left(\dfrac{X_1+Y_1\sqrt{d}}{2}\right)^k\qquad {\text{\rm for~all}} \quad k\ge 1.\label{Pellsoln2}
\end{eqnarray*}
Furthermore, the sequences $ \{x_k\}_{k\ge 1} $ and $ \{X_k\}_{k\geq 1} $ are binary recurrent. In fact, the following formulae
\begin{eqnarray*}
x_k=\dfrac{(x_1+y_1\sqrt{d})^{k}+(x_1-y_1\sqrt{d})^{k}}{2},
\end{eqnarray*}
and 
\begin{eqnarray*}
X_k = \left(\dfrac{X_1+Y_1\sqrt{d}}{2}\right)^{k}+ \left(\dfrac{X_1-Y_1\sqrt{d}}{2}\right)^{k}
\end{eqnarray*}
hold for all positive integers $ k $.

Recently, Bravo, G\'omez-Ruiz and Luca \cite{Bravo1} studied the Diophantine equation
\begin{eqnarray}\label{Bravo}
x_l = T_m+T_n,
\end{eqnarray}
where $ x_l $ are the $ x- $coordinates of the solutions of the Pell equation \eqref{Pelleqn1} for some positive integer $ l $ and $ \{T_n\}_{n\geq 0} $ is the sequence of Tribonacci numbers given by $ T_0=0,~ T_1=1=T_2 $ and $ T_{n+3} = T_{n+2}+T_{n+1}+T_n $ for all $ n\ge 0 $. They proved that for each square free integer $ d\ge 2 $, there is at most one positive integer $ l $ such that $ x_l $ admits the representation \eqref{Bravo} for some nonnegative integers $ 0\le m\le n $, except for $ d\in\{2,3,5,15,26\} $. Furthermore, they explicitly stated all the solutions for these exceptional cases. 

In the same spirit, Rihane, Hernane and Togb\'e \cite{Togbe} studied the Diophantine equations
\begin{eqnarray}\label{Togbe}
x_n = P_m \qquad \text{and}\qquad X_n=P_m,
\end{eqnarray}
where $ x_n $ and $ X_n $ are the $ x- $coordinates of the solutions of the Pell equations \eqref{Pelleqn1} and \eqref{Pelleqn2}, respectively, for some positive integers $ n $ and  $ \{P_m\}_{m\ge 0} $ is the sequence of Padovan numbers. They proved that for each square free integer $ d\ge 2 $, there is at most one positive integer $ x $  participating in the Pell equation \eqref{Pelleqn1} and one positive integer $ X $ participating in the Pell equation \eqref{Pelleqn2} that is a Padovan number with a few exceptions of $ d $ that they effectively computed. Furthermore, the exceptional cases were $ d\in\{2,3,5,6 $ and $ d\in\{5\} $  for the the first  and second equations  in \eqref{Togbe}, respectively. Several other related problems have been studied where $ x_l $ belongs to some interesting positive integer sequences. For example, see \cite{DM1, Dossavi, Faye, BLT1, BLT2, BLT3, Luca16, Luca15}.
\section{Main Results}
In this paper, we study a problem related to that of Bravo, G\'omez-Ruiz and Luca \cite{Bravo1} but with the Padovan sequence instead of the Tribonacci sequence. We also extend the results from the Pell equation \eqref{Pelleqn1} to the Pell equation \eqref{Pelleqn2}. In both cases we find that there are only finitely many solutions that we  effectively compute.

Since $ P_1=P_2=P_3 = 1$, we discard the situations when $ n=1 $ and $ n=2 $ and just count the solutions for $ n=3 $. Similarly, $ P_4=P_5=2 $, we discard the situation when $ n=4 $ and just count the solutions for $ n=5 $. The main aim of this paper is to prove the following results.

\begin{theorem}\label{Main1} 
For each integer $ d\ge 2 $ which is not a square, there is atmost one positive integer $ k $ such that $ x_k $ admits a representation as
\begin{eqnarray}
x_k=P_n+P_m
\end{eqnarray}
for some nonnegative integers $ 0\le m\le n $, except when $ d\in\{2,3,6,15,110,483\} $ in the $ +1 $ case and $ d\in\{2,5,10,17\} $ in the $ -1 $ case.
\end{theorem}

\begin{theorem}\label{Main2}
For each integer $ d\ge 2 $ which is not a square, there is atmost one positive integer $ k $ such that $ X_k $ admits a representation as
\begin{eqnarray}
X_k=P_n+P_m
\end{eqnarray}
for some nonnegative integers $ 0 \le m\le n $, except when $ d\in\{3,5,21\}$ in the $ +4 $ case and $ d\in\{2, 5\} $ in the $ -4 $ case.
\end{theorem}
For the exceptional values of $ d $ listed in Theorem \ref{Main1} and Theorem \ref{Main2}, all solutions $ (k,n,m) $ are listed at the end of the proof of each result. The main tools used in this paper are the lower bounds for linear forms in logarithms of algebraic numbers and the Baker-Davenport reduction procedure, as well as the elementary properties of Padovan numbers and solutions to Pell equations.

\section{Preliminary results}
\subsection{The Padovan sequence} 
Here, we recall some important properties of the Padovan sequence $ \{P_n\}_{n\geq 0} $. The characteristic equation
\begin{eqnarray*}
 \Psi (x):= x^3-x-1 = 0
\end{eqnarray*}
has roots $ \alpha, \beta, \gamma = \bar{\beta} $, where
\begin{eqnarray}\label{Pado1}
\alpha =\dfrac{r_1+r_2}{6}, \qquad \beta = \dfrac{-(r_1+r_2)+\sqrt{-3}(r_1-r_2)}{12} 
\end{eqnarray}
and
\begin{eqnarray}\label{Pado2}
r_1=\sqrt[3]{108+12\sqrt{69}} \quad \text{and}\quad r_2=\sqrt[3]{108-12\sqrt{69}}.
\end{eqnarray}
Furthermore, the Binet formula is given by
\begin{eqnarray}\label{Pado3}
P_n = a\alpha^{n}+b\beta^{n}+c\gamma^{n} \qquad \text{ for all} \quad n\ge 0,
\end{eqnarray}
where
\begin{eqnarray}\label{Pado4}
a=\dfrac{(1-\beta)(1-\gamma)}{(\alpha-\beta)(\alpha-\gamma)}, \quad b= \dfrac{(1-\alpha)(1-\gamma)}{(\beta -\alpha)(\beta-\gamma)}, \quad c = \dfrac{(1-\alpha)(1-\beta)}{(\gamma-\alpha)(\gamma-\beta)}=\bar{b}.
\end{eqnarray}
Numerically, the following estimates hold:
\begin{eqnarray}\label{Pado5}
&1.32<\alpha<1.33\nonumber\\
&0.86 < |\beta|=|\gamma|=\alpha^{-\frac{1}{2}}< 0.87\\
&0.72<a<0.73\nonumber\\
&0.24<|b|=|c|<0.25.\nonumber
\end{eqnarray}
From \eqref{Pado1}, \eqref{Pado2} and \eqref{Pado5}, it is easy to see that the contribution the complex conjugate roots $ \beta $ and $ \gamma $, to the right-hand side of \eqref{Pado3}, is very small. In particular, setting
\begin{eqnarray}\label{Pado6}
e(n):=P_n-a\alpha^{n}=b\beta^{n}+c\gamma^{n}\quad \text{ then } \quad |e(n)|< \dfrac{1}{\alpha^{n/2}}
\end{eqnarray}
holds for all $ n\ge 1 $.
Furthermore, by induction, we can prove that 
\begin{eqnarray}\label{Pado7}
\alpha^{n-2}\leq P_n \leq \alpha^{n-1} \quad \text{holds for all }\quad n\geq 4.
\end{eqnarray}

\subsection{Linear forms in logarithms}
Let $ \eta $ be an algebraic number of degree $ D $ with minimal primitive polynomial over the integers
$$ a_{0}x^{d}+ a_{1}x^{d-1}+\cdots+a_{d} = a_{0}\prod_{i=1}^{D}(x-\eta^{(i)}),$$
where the leading coefficient $ a_{0} $ is positive and the $ \eta^{(i)} $'s are the conjugates of $ \eta $. Then the \textit{logarithmic height} of $ \eta $ is given by
$$ h(\eta) := \dfrac{1}{D}\left( \log a_{0} + \sum_{i=1}^{d}\log\left(\max\{|\eta^{(i)}|, 1\}\right)\right).$$

In particular, if $ \eta = p/q $ is a rational number with $ \gcd (p,q) = 1 $ and $ q>0 $, then $ h(\eta) = \log\max\{|p|, q\} $. The following are some of the properties of the logarithmic height function $ h(\cdot) $, which will be used in the next sections of this paper without reference:
\begin{eqnarray}
h(\eta_1\pm \eta_2) &\leq& h(\eta) +h(\eta_1) +\log 2,\nonumber\\
h(\eta_1\eta_2^{\pm 1})&\leq & h(\eta_1) + h(\eta_2),\\
h(\eta^{s}) &=& |s|h(\eta) ~~~~~~ (s\in\mathbb{Z}). \nonumber
\end{eqnarray}

\begin{theorem}\label{Matveev11} Let $\eta_1,\ldots,\eta_t$ be positive real algebraic numbers in a real algebraic number field  $\mathbb{K} \subset \mathbb{R}$ of degree $D_\mathbb{K}$, $b_1,\ldots,b_t$ be nonzero integers, and assume that
\begin{equation}
\label{eq:Lambda}
\Lambda:=\eta_1^{b_1}\cdots\eta_t^{b_t} - 1,
\end{equation}
is nonzero. Then
$$
\log |\Lambda| > -1.4\times 30^{t+3}\times t^{4.5}\times D_\mathbb{K}^{2}(1+\log D_{\mathbb{K}})(1+\log B)A_1\cdots A_t,
$$
where
$$
B\geq\max\{|b_1|, \ldots, |b_t|\},
$$
and
$$A
_i \geq \max\{D_{\mathbb{K}} h(\eta_i), |\log\eta_i|, 0.16\},\qquad {\text{for all}}\qquad i=1,\ldots,t.
$$
\end{theorem} 
\subsection{Reduction procedure}\label{Reduction}
During the calculations, we get upper bounds on our variables which are too large, thus we need to reduce them. To do so, we use some results from the theory of continued fractions. 

For the treatment of linear forms homogeneous in two integer variables, we use the well-known classical result in the theory of Diophantine approximation.
\begin{lemma}\label{Legendre}
Let $\tau$ be an irrational number,  $ \frac{p_0}{q_0}, \frac{p_1}{q_1}, \frac{p_2}{q_2}, \ldots $ be all the convergents of the continued fraction of $ \tau $ and $ M $ be a positive integer. Let $ N $ be a nonnegative integer such that $ q_N> M $. Then putting $ a(M):=\max\{a_{i}: i=0, 1, 2, \ldots, N\} $, the inequality
\begin{eqnarray*}
\left|\tau - \dfrac{r}{s}\right|> \dfrac{1}{(a(M)+2)s^{2}},
\end{eqnarray*}
holds for all pairs $ (r,s) $ of positive integers with $ 0<s<M $.
\end{lemma}

For a nonhomogeneous linear form in two integer variables, we use a slight variation of a result due to Dujella and Peth{\H o} (see \cite{dujella98}, Lemma 5a). For a real number $X$, we write  $||X||:= \min\{|X-n|: n\in\mathbb{Z}\}$ for the distance from $X$ to the nearest integer.
\begin{lemma}\label{Dujjella}
Let $M$ be a positive integer, $\frac{p}{q}$ be a convergent of the continued fraction of the irrational number $\tau$ such that $q>6M$, and  $A,B,\mu$ be some real numbers with $A>0$ and $B>1$. Let further 
$\varepsilon: = ||\mu q||-M||\tau q||$. If $ \varepsilon > 0 $, then there is no solution to the inequality
$$
0<|u\tau-v+\mu|<AB^{-w},
$$
in positive integers $u,v$ and $w$ with
$$ 
u\le M \quad {\text{and}}\quad w\ge \dfrac{\log(Aq/\varepsilon)}{\log B}.
$$
\end{lemma}

At various occasions, we need to find a lower bound for linear forms in logarithms with bounded integer coefficients in three and four variables. In this case we use the LLL- algorithm that we describe below. Let $ \tau_1, \tau_2, \ldots \tau_t \in\mathbb{R}$ and the linear form
\begin{eqnarray}
x_1\tau_1+x_2\tau_2+\cdots+x_t\tau_t \quad \text{ with } \quad |x_i|\le X_i.
\end{eqnarray}
We put $ X:=\max\{X_i\} $, $ C> (tX)^{t} $ and consider the integer lattice $ \Omega $ generated by
\begin{eqnarray*}
\textbf{b}_j: = \textbf{e}_j+\lfloor C\tau_j\rceil \quad \text{ for} \quad 1\le j\le t-1 \quad \text{ and} \quad \textbf{b}_t:=\lfloor C\tau_t\rceil \textbf{e}_t,
\end{eqnarray*}
where $ C $ is a sufficiently large positive constant.
\begin{lemma}\label{LLL}
Let $ X_1, X_2, \ldots, X_t $ be positive integers such that $ X:=\max\{X_i\} $ and $ C> (tX)^{t} $ is a fixed sufficiently large constant. With the above notation on the lattice $ \Omega $, we consider a reduced base $ \{\textbf{b}_i \}$ to $ \Omega $ and its associated Gram-Schmidt orthogonalization base $ \{\textbf{b}_i^*\}$. We set
\begin{eqnarray*}
c_1:=\max_{1\le i\le t}\dfrac{||\textbf{b}_1||}{||\textbf{b}_i^*||}, \quad \theta:=\dfrac{||\textbf{b}_1||}{c_1}, \quad Q:=\sum_{i=1}^{t-1}X_i^{2} \quad \text{and} \quad R:=\left(1+ \sum_{i=1}^{t}X_i\right)/2.
\end{eqnarray*}
If the integers $ x_i $ are such that $ |x_i|\le X_i $, for $ 1\le i \le t $ and $ \theta^2\ge Q+R^2 $, then we have
\begin{eqnarray*}
\left|\sum_{i=1}^{t}x_i\tau_i\right|\ge \dfrac{\sqrt{\theta^2-Q}-R}{C}.
\end{eqnarray*}
\end{lemma}
\noindent
For the proof and further details, we refer the reader to the book of Cohen. (Proposition 2.3.20 in [\cite{Cohen}, Pg. 58--63).

Finally, the following Lemma is also useful. It is Lemma 7 in \cite{guzmanluca}. 
\begin{lemma}
\label{gl}
If $r\geqslant 1$, $H>(4r^2)^r$  and $H>L/(\log L)^r$, then
$$
L<2^rH(\log H)^r.
$$
\end{lemma}


%

\section{Proof of Theorem \ref{Main1}}
Let $ (x_1, y_1) $ be the smallest positive integer solution to the Pell quation \eqref{Pelleqn1}. We Put
\begin{eqnarray}\label{Kay1}
\delta:=x_1+y_1\sqrt{d} \quad \text{and} \quad \sigma=x_1-y_1\sqrt{d}.
\end{eqnarray}
From which we get that
\begin{eqnarray}\label{kalayi1}
\delta\cdot\sigma=x_1^2-dy_1^2 =: \epsilon, \quad \text{where} \quad \epsilon\in\{\pm 1\}. 
\end{eqnarray}
Then
\begin{eqnarray}\label{Kayii1}
x_k = \dfrac{1}{2}(\delta^k+\sigma^k).
\end{eqnarray}
Since $ \delta \ge 1+\sqrt{2} $, it follows that the estimate 
\begin{eqnarray}\label{kalayim1}
\dfrac{\delta^{k}}{\alpha^{4}}\le x_k \le \delta^k \quad \text{ holds for all } \quad k\ge 1.
\end{eqnarray}
We assume that $ (k_1, n_1, m_1) $ and $ (k_2, n_2, m_2) $ are triples of integers such that
\begin{eqnarray}\label{kalayim2}
x_{k_1}=P_{n_1}+P_{m_1} \quad \text{and} \quad x_{k_2}=P_{n_2}+P_{m_2}
\end{eqnarray}
We asuume that $ 1\le k_1 < k_2 $. We also assume that $ 3\le m_i< n_i $ for $ i=1,2 $.
We set $ (k,n,m):=(k_i, n_i, m_i) $, for $ i=1,2 $. Using the inequalities \eqref{Pado7} and \eqref{kalayim1}, we get from \eqref{kalayim2} that 
\begin{eqnarray*}
\dfrac{\delta^{k}}{\alpha^{4}}\le x_k=P_n+P_m\le 2\alpha^{n-1}\quad \text{and} \quad \alpha^{n-2}\le P_n+P_m = x_k\le \delta^{k}.
\end{eqnarray*}
The above inequalities give
\begin{eqnarray*}
(n-2)\log\alpha<k\log\delta < (n+3)\log\alpha +\log 2.
\end{eqnarray*}
Dividing through by $ \log\alpha $ and setting $ c_2:=1/\log\alpha $, we get that 
\begin{eqnarray*}
-2<c_2k\log\delta-n<3+c_2\log2,
\end{eqnarray*}
and since $ \alpha^3>2 $, we get
\begin{eqnarray}\label{kabizi}
|n-c_2k\log\delta|<6.
\end{eqnarray}
Furthermore, $ k<n $, for if not, we would then get that 
\begin{eqnarray*}
\delta^{n}\le \delta^{k}<2\alpha^{n+3}, \quad \text{implying}\quad \left(\dfrac{\delta}{\alpha}\right)^{n}<2\alpha^{3},
\end{eqnarray*}
which is false since $ \delta\ge 1+\sqrt{2} $, $ 1.32<\alpha<1.33 $ (by \eqref{Pado5}) and $ n\ge 4 $.

Besides, given that $ k_1<k_2 $, we have by \eqref{Pado7} and \eqref{kalayim2} that
\begin{eqnarray*}
\alpha^{n_1-2}\le P_{n_1}\le P_{n_1}+P_{m_1}=x_{k_1}<x_{k_2}= P_{n_2}+P_{m_2} \le 2P_{n_2}<2\alpha^{n_2-1}.
\end{eqnarray*}
Thus, we get that
\begin{eqnarray}\label{faji1}
n_1<n_2+4.
\end{eqnarray}
\subsection{An inequality for $ n $ and $ k $ (I)}
Using the equations \eqref{Pado3} and \eqref{Kayii1} and \eqref{kalayim2}, we get
\begin{eqnarray*}
\dfrac{1}{2}(\delta^{k}+\sigma^{k})=P_n+P_m=a\alpha^{n}+e(n)+a\alpha^{m}+e(m)
\end{eqnarray*}
So,
\begin{eqnarray*}
\frac{1}{2}\delta^{k}-a(\alpha^{n}+\alpha^{m})=-\dfrac{1}{2}\sigma^{k}+e(n)+e(m),
\end{eqnarray*}
and by \eqref{Pado6}, we have
\begin{eqnarray*}
\left|\delta^{k}(2a)^{-1}\alpha^{-n}(1+\alpha^{m-n})^{-1}-1\right|&\leq& \dfrac{1}{2\delta^{k}a(\alpha^{n}+\alpha^{m})}+\dfrac{2|b|}{\alpha^{n/2}a(\alpha^{n}+\alpha^{m})}\\&&+\dfrac{2|b|}{\alpha^{m/2}a(\alpha^{n}+\alpha^{m})}\\
&\le& \dfrac{1}{a\alpha^{n}}\left(\dfrac{1}{2\delta^{k}}+\dfrac{2|b|}{\alpha^{n/2}}+\dfrac{2|b|}{\alpha^{m/2}}\right)<\dfrac{1.5}{\alpha^{n}}.
\end{eqnarray*}
Thus, we have
\begin{eqnarray}\label{kapuli1}
\left|\delta^{k}(2a)^{-1}\alpha^{-n}(1+\alpha^{m-n})^{-1}-1\right|&<&\dfrac{1.5}{\alpha^{n}}.
\end{eqnarray}
Put
\begin{eqnarray*}
\Lambda_1:=\delta^{k}(2a)^{-1}\alpha^{-n}(1+\alpha^{m-n})^{-1}-1
\end{eqnarray*}
and 
\begin{eqnarray*}
\Gamma_1:=k\log\delta-\log(2a) -n\log\alpha-\log(1+\alpha^{m-n}).
\end{eqnarray*}
Since $|\Lambda_1|=|e^{\Gamma_1}-1|<\frac{1}{2}$ for $ n\ge 4 $ (because $1.5/\alpha^{4}<1/2$), since the inequality $|y|<2|e^{y}-1|$ holds for all $ y\in \left(-\frac{1}{2}, \frac{1}{2}\right) $, it follows that $ e^{|\Gamma_1|}<2 $ and so
\begin{eqnarray*}
|\Gamma_1|<e^{|\Gamma_1|}|e^{\Gamma_1}-1|<\dfrac{3}{\alpha^{n}}.
\end{eqnarray*}
Thus, we get that
\begin{eqnarray}\label{kabuzi}
\left|k\log\delta-\log(2a) -n\log\alpha-\log(1+\alpha^{m-n})\right|<\dfrac{3}{\alpha^{n}}.
\end{eqnarray}
We apply Theorem \ref{Matveev11} on the left-hand side of \eqref{kapuli1} with the data:
\begin{eqnarray*}
&t:=4, \quad \eta_{1}:=\delta, \quad \eta_2:=2a, \quad \eta_3:=\alpha, \eta_4: =1+\alpha^{m-n},\\
&b_1:=k, \quad b_2:=-1, \quad b_3:=-n, \quad b_4:=-1.
\end{eqnarray*}
Furthermore, we take the number field $ \mathbb{K}=\mathbb{Q}(\sqrt{d}, \alpha) $ which has degree $ D=6 $. Since $ \max\{1,k,n\}\leq n $, we take $D_{\mathbb{K}}=n$.
First we note that the left-hand side of \eqref{kapuli1} is non-zero, since otherwise,
\begin{eqnarray*}
\delta^{k}=2a(\alpha^{n}+\alpha^{m}).
\end{eqnarray*}
The left-hand side belongs to the quadratic field $ \mathbb{Q}(\sqrt{d}) $ while the right-hand side belongs to the cubic field $ \mathbb{Q}(\alpha) $. These fields only intersect when both sides are rational numbers. Since $ \delta^{k} $ is a positive algebraic integer and a unit, we get that  to $ \delta^{k} =1 $. Hence, $ k=0 $, which is a contradiction. Thus, $ \Lambda_1\neq 0 $ and we can apply Theorem \ref{Matveev11}.

We have $ h(\eta_1)=h(\delta)=\frac{1}{2}\log\delta $ and $ h(\eta_3)=h(\alpha)=\frac{1}{3}\log\alpha $. Further, $$ 2a=\dfrac{2\alpha(\alpha+1)}{3\alpha^2-1}, $$ the mimimal polynomial of $ 2a $ is $ 23x^3-46x^2+24x-8 $ and has roots $ 2a, 2b, 2c $. Since $ 2|b|=2|c|<1 $ (by \eqref{Pado5}), then
\begin{eqnarray*}
h(\eta_2)=h(2a)=\dfrac{1}{3}(\log 23+\log (2a)).
\end{eqnarray*}
On the other hand,
\begin{eqnarray*}
h(\eta_4)&=&h(1+\alpha^{m-n})\leq h(1)+h(\alpha^{m-n})+\log 2\\
&=&(n-m)h(\alpha)+\log2 = \dfrac{1}{3}(n-m)\log\alpha +\log2.
\end{eqnarray*}
Thus, we can take $ A_1:=3\log\delta $,
\begin{eqnarray*}
 A_2:=2(\log 23+\log (2a)), \quad A_3:=2\log\alpha, \quad A_4:=2(n-m)\log\alpha+6\log 2.
\end{eqnarray*}
Now, Theorem \ref{Matveev11} tells us that
\begin{eqnarray*}
\log|\Lambda_1|&> &-1.4\times 30^{7}\times 4^{4.5}\times 6^{2}(1+\log 6)(1+\log n)(3\log\delta)\\
&& \times (2(\log 23+\log(2a))(2\log\alpha)(2(n-m)\log\alpha+6\log 2)\\
&>&-2.33\times 10^{17}(n-m)(\log n)(\log\delta).
\end{eqnarray*}
Comparing the above inequality with \eqref{kapuli1}, we get
\begin{eqnarray*}
n\log\alpha - \log 1.5 < 2.33\times 10^{17}(n-m)(\log n)(\log\delta).
\end{eqnarray*}
Hence, we get that
\begin{eqnarray}\label{good1}
n<8.30\times 10^{17}(n-m)(\log n)(\log\delta).
\end{eqnarray}
We now return to the equation $ x_k=P_n+P_m $ and rewrite it as
\begin{eqnarray*}
\dfrac{1}{2}\delta^{k}-a\alpha^{n} = -\dfrac{1}{2}\sigma^{k}+e(n)+P_m,
\end{eqnarray*}
we obtain
\begin{eqnarray}\label{kayija1}
\left|\delta^{k}(2a)^{-1}\alpha^{-n}-1\right|\leq \dfrac{1}{a\alpha^{n-m}}\left(\dfrac{1}{\alpha}+\dfrac{1}{\alpha^{m+n/2}}+\dfrac{1}{2\delta^{k}\alpha^{m}}\right)<\dfrac{2.5}{\alpha^{n-m}}.
\end{eqnarray}
Put
\begin{eqnarray*}
\Lambda_2:=\delta^{k}(2a)^{-1}\alpha^{-n}-1, \quad \Gamma_{2}:=k\log\delta-\log(2a)-n\log\alpha.
\end{eqnarray*}
We assume for technical reasons that $ n-m\ge 10 $. So $ |e^{\Lambda_2}-1|<\frac{1}{2} $. It follows that 
\begin{eqnarray}\label{mukazi1}
\left|k\log\delta-\log(2a)-n\log\alpha\right|=|\Gamma_2|<e^{|\Lambda_2|}|e^{\Lambda_2}-1|<\dfrac{5}{\alpha^{n-m}}.
\end{eqnarray}
Furthermore, $ \Lambda_2\neq 0 $ (so $\Gamma_2 \neq 0$), since $ \delta^{k}\in\mathbb{Q}(\alpha) $ by the previous argument.

We now apply Theorem \ref{Matveev11} to the left-hand side of \eqref{kayija1} with the data
\begin{eqnarray*}
t:=3, \quad \eta_1:=\delta, \quad\eta_2:=2a, \quad \eta_3:=\alpha, \quad b_1:=k, \quad b_2:=-1, \quad b_3:=-n.
\end{eqnarray*}
Thus, we have the same $ A_1, ~A_2, A_3 $ as before. Then, by Theorem \ref{Matveev11}, we conclude that
\begin{eqnarray*}
\log|\Lambda|>-9.82\times 10^{14}(\log\delta)(\log n)(\log\alpha).
\end{eqnarray*}
By comparing with \eqref{kayija1}, we get
\begin{eqnarray}\label{good2}
n-m <9.84\times 10^{14}(\log\delta)(\log n).
\end{eqnarray}
This was obtained under the assumption that $ n-m\ge 10 $, but if $ n-m<10$, then the inequality also holds as well. We replace the bound \eqref{good2} on $ n-m $ in \eqref{good1} and use the fact that $ \delta^{k}\le 2\alpha^{n+3} $, to obtain bounds on $ n $ and $ k $ in terms of $ \log n $ and $\log\delta$. We now record what we have proved so far.
\begin{lemma}
Let $ (k,n,m) $ be a solution to the equation  $x_k=P_n+P_m$ with $ 3\le m<n $, then
\begin{eqnarray}\label{lemmata1}
k< 2.5\times 10^{32}(\log n)^{2}(\log\delta) \quad \text{and} \quad n<8.2\times 10^{32}(\log n)^{2}(\log\delta)^{2}.
\end{eqnarray}
\end{lemma}
\subsection{Absolute bounds (I)}
We recall that $ (k,n,m)=(k_i,n_i, m_i) $, where $ 3\le m_i<n_i $, for $ i=1,2 $ and $ 1\le k_1<k_2 $. Further, $ n_i\ge 4 $ for $ i=1,2 $. We return to \eqref{mukazi1} and write
\begin{eqnarray*}
\left|\Gamma_2^{(i)}\right|:=\left|k_i\log\delta - \log(2a) -n_i\log\alpha\right|<\dfrac{5}{\alpha^{n_i-m_i}}, \quad \text{ for } \quad i=1,2.
\end{eqnarray*}
We do a suitable cross product between $ \Gamma_2^{(1)}, ~ \Gamma_2^{(2)} $ and $ k_1, k_2 $ to eliminate the term involving $ \log\delta $ in the above linear forms in logarithms:
\begin{eqnarray}
|\Gamma_3|&:=&|(k_1-k_2)\log(2a)+(k_1n_2-k_2n_1)\log\alpha|=|k_2\Gamma_2^{(1)}-k_1\Gamma_2^{(2)}|\nonumber\\
&\le& k_2|\Gamma_2^{(1)}|+k_1|\Gamma_2^{(2)}|\quad \le \quad  \dfrac{5k_2}{\alpha^{n_1-m_1}}+\dfrac{5k_1}{\alpha^{n_2-m_2}}\quad\le \quad \dfrac{10n_2}{\alpha^{\lambda}},\label{basaja1}
\end{eqnarray}
where \[ \lambda:=\min_{1\le i\leq 2} \{n_i-m_i\}.\]

We need to find an upper bound for $ \lambda $. If $ 10n_2/\alpha^{\lambda} > 1/2 $, we then get
\begin{eqnarray}\label{kabaya122}
\lambda < \dfrac{\log (20n_2)}{\log \alpha}<4\log(20n_2).
\end{eqnarray}
Otherwise, $ |\Gamma_3|<\frac{1}{2} $, so
\begin{eqnarray}\label{kabaya22}
\left|e^{\Gamma_3}-1\right|=\left|(2a)^{k_1-k_2}\alpha^{k_1n_2-k_2n_1}-1\right|<2|\Gamma_3|<\dfrac{20n_2}{\alpha^{\lambda}}.
\end{eqnarray}
We apply Theorem \ref{Matveev11} with the data: $ t:=2 $, $ \eta_1 := 2a$, $ \eta_2:= \alpha$, $ b_1:=k_1-k_2 $, $ b_2:=k_1n_2-k_2n_1 $. We take the number field $ \mathbb{K}:=\mathbb{Q}(\alpha) $ and $ D=3 $. We begin by checking that $ e^{\Gamma_3}-1\neq 0 $ (so $ \Gamma_3\neq 0 $). This is true because $ \alpha $ and $ 2a $ are multiplicatively independent, since $ \alpha $ is a unit in the ring of integers $ \mathbb{Q}(\alpha) $ while the norm of $ 2a $ is $ 8/23 $.

We note that $ |k_1-k_2|<k_2<n_2 $. Further, from \eqref{basaja1}, we have
\begin{eqnarray*}
|k_2n_1-k_1n_2|<(k_2-k_1)\dfrac{|\log(2a)|}{\log\alpha}+\dfrac{10k_2}{\alpha^{\lambda}\log\alpha}<11k_2<11n_2
\end{eqnarray*}
given that $ \lambda\ge 1 $. So, we can take $ B:=11n_2 $.
By Theorem \ref{Matveev11}, with the same $ A_1:=\log 23 $ and $ A_2:=\log\alpha $, we have that 
\begin{eqnarray*}
\log|e^{\Gamma_3}-1|>-1.55\times 10^{11}(\log n_2)(\log\alpha).
\end{eqnarray*}
By comparing this with \eqref{kabaya22}, we get
\begin{eqnarray}\label{kabaya123}
\lambda <1.56\times 10^{11}\log n_2.
\end{eqnarray}
Note that \eqref{kabaya123} is better than \eqref{kabaya122}, so \eqref{kabaya123} always holds. Without loss of generality, we can assume that $ \lambda = n_i-m_i $, for $ i=1,2 $ fixed.

We set $ \{i,j\}=\{1,2\} $ and return to \eqref{kabuzi} to replace $ (k,n,m)=(k_i,n_i, m_i) $:
\begin{eqnarray}\label{muka11}
|\Gamma_1^{(i)}|=\left|k_i\log\delta-\log(2a) -n_i\log\alpha-\log(1+\alpha^{m_i-n_i})\right|<\dfrac{3}{\alpha^{n_i}},
\end{eqnarray}
and also return to \eqref{mukazi1}, replacing with $ (k, n,m)=(k_j, n_j, m_j) $:
\begin{eqnarray}\label{muka12}
|\Gamma_2^{(j)}|=\left|k_j\log\delta-\log(2a)-n_j\log\alpha\right|<\dfrac{5}{\alpha^{n_j-m_j}}.
\end{eqnarray}
We perform a cross product on \eqref{muka11} and \eqref{muka12} in order to eliminate the term on $ \log\delta $:
\begin{eqnarray}
|\Gamma_4|&:=&\left|(k_j-k_i)\log(2a)+(k_jn_i-k_in_j)\log\alpha+k_j\log(1+\alpha^{m_i-n_i})\right|\nonumber\\
&=&\left|k_i\Gamma_2^{(j)}-k_j\Gamma_1^{(i)}\right|\le k_i\left|\Gamma_2^{(j)}\right|+k_j\left|\Gamma_1^{(i)}\right|\nonumber\\
&<&\dfrac{5k_i}{\alpha^{n_j-m_j}}+\dfrac{3k_j}{\alpha^{n_i}}<\dfrac{8n_2}{\alpha^{\nu}}\label{kipro1}
\end{eqnarray}
with $ \nu:=\min\{n_i, n_j-m_j\} $. As before, we need to find an upper bound on $ \nu $. If $ 8n_2/\alpha^{\nu}>1/2 $, then we get
\begin{eqnarray}\label{boss1}
\nu < \dfrac{\log (16n_2)}{\log\alpha}< 4\log (16n_2).
\end{eqnarray}
Otherwise, $ |\Gamma_4|<1/2 $, so we have
\begin{eqnarray}\label{bosco1}
\left|e^{\Gamma_4}-1\right|\le 2|\Gamma_4|<\dfrac{16n_2}{\alpha^{\nu}}.
\end{eqnarray}
In order to apply Theorem \ref{Matveev11}, first if $ e^{\Gamma_4}=1 $, we obtain
\begin{eqnarray*}
(2a)^{k_i-k_j}=\alpha^{k_jn_i-k_in_j}(1+\alpha^{-\lambda})^{k_j}.
\end{eqnarray*}
Since $ \alpha $ is a unit, the right-hand side in above is an algebraic integer. This is a contradiction because $ k_1<k_2 $ so $ k_i-k_j\neq 0 $, and neither $ (2a) $ nor $ (2a)^{-1} $ are algebraic intgers. Hence $ e^{\Gamma_4}\neq 1 $. By assuming that $ \nu \ge 100 $, we apply Theorem \ref{Matveev11} with the data:
\begin{eqnarray*}
&t:=3, \quad \eta_1:=2a, \quad \eta_2:=\alpha, \quad \eta_3:=1+\alpha^{-\lambda},\\
& b_1:=k_j-k_i, \quad b_2:=k_jn_i-k_in_j, \quad b_3:=k_j,
\end{eqnarray*} 
and the inequalities \eqref{kabaya123} and \eqref{bosco1}. We get
\begin{eqnarray}
\nu=\min\{n_i, n_j-m_j\}<1.14\times 10^{14}\lambda\log n_2<1.78\times 10^{25}(\log n_2)^{2}.
\end{eqnarray}
The above inequality also holds when $ \nu <100 $. Further, it also holds when the inequality \eqref{boss1} holds. So the above inequality holds in all cases. Note that the case $ \{i, j\} =\{2, 1\} $ leads to $ n_1-m_1\le n_1\le n_2+4 $ whereas $ \{i, j\} = \{1,2\} $ lead to $ \nu = \min\{n_1, n_2-m_2\} $. Hence, either the minimum is $  n_1 $, so
\begin{eqnarray}\label{case1}
n_1< 1.78\times 10^{25}(\log n_2)^{2},
\end{eqnarray}  
or the minimum is $ n_j-m_j $ and from the inequality \eqref{kabaya123} we get that
\begin{eqnarray}\label{case2}
\max_{1\le j\le 2}\{n_j- m_j\}< 1.78\times 10^{25}(\log n_2)^{2}.
\end{eqnarray}
Next, we assume that we are in the case \eqref{case2}. We evaluate \eqref{muka11} in $ i=1,2 $ and make a suitable cross product to eliminate the term involving $ \log\delta $:
\begin{eqnarray}
\left|\Gamma_5\right|&:=&\left|(k_2-k_1)\log(2a)+(k_2n_1-k_1n_2)\log\alpha \right.\nonumber\\&&\left.+k_2\log (1+\alpha^{m_1-n_1})-k_1\log (1+\alpha^{m_2-n_2})\right|\nonumber\\
&=&\left|k_1\Gamma_1^{(2)}-k_2\Gamma_1^{(1)}\right|\le k_1\left|\Gamma_1^{(2)}\right|+k_2\left|\Gamma_1^{(1)}\right|<\dfrac{6n_2}{\alpha^{n_1}}.\label{kakawu1}
\end{eqnarray}
In the above inequality we used the inequality \eqref{faji1}to conclude that $ \min\{n_1, n_2\}\ge n_1-4 $ as well as the fact that $ n_i\ge 4 $ for $ i=1.2 $. Next, we apply a linear form in four logarithms to obtain an upper bound to $ n_1 $. As in the previous calculations, we pass from \eqref{kakawu1} to
\begin{eqnarray}\label{kakawu2}
\left|e^{\Gamma_5}-1\right|<\dfrac{12n_2}{\alpha^{n_1}},
\end{eqnarray}
which is implied by \eqref{kakawu1} except if $ n_1 $ is very small, say
\begin{eqnarray}\label{kakawu3}
n_1\le 4\log(12n_2).
\end{eqnarray}
Thus, we assume that \eqref{kakawu3} does not hold, therefore \eqref{kakawu2} holds. Then to apply Theorem \ref{Matveev11}, we fist justify that $ e^{\Gamma_5}\neq 1 $. Otherwise, 
\begin{eqnarray*}
(2a)^{k_1-k_2}=\alpha^{k_2n_1-k_1n_2}(1+\alpha^{n_1-m_1})^{k_2}(1+\alpha^{n_2-m_2})^{-k_1},
\end{eqnarray*}
By the fact that $ k_1<k_2 $, the norm  $ \textbf{N}_{\mathbb{Q}(\alpha)/\mathbb{Q}} (2a) =\frac{8}{23} $ and that $ \alpha $ is a unit, we have that $ 23 $ divides the norm  $ \textbf{N}_{\mathbb{K}/\mathbb{Q}} (1+\alpha^{n_1-m_1}) $. The factorization of the ideal generated by $ 23 $ in $ \mathcal{O}_{\mathbb{Q}(\alpha)} $ is $ (23)=\mathfrak{p}_1^2\mathfrak{p}_2 $, where $ \mathfrak{p}_1=(23, ~\alpha+13) $ and $ \mathfrak{p}_2=(23,~ \alpha+20) $. Hence $ \mathfrak{p}_2 $ divides $ \alpha^{n_1-m_1} + 1 $. Given that $ \alpha \equiv -20~ (\text{mod}~ \mathfrak{p}_2) $, then $ (-20)^{n_1-m_1}\equiv -1(\text{mod}~\mathfrak{p}_2)  $. Taking the norm $ \textbf{N}_{\mathbb{Q}(\alpha)/\mathbb{Q}} $, we obtain that $ (-20)^{n_1-m_1} \equiv -1~(\text{mod}~ 23) $. If $ n_1-m_1 $ is even $ -1 $ is a quadratic residue modulo $ 23 $ and if $ n_1-m_1 $ is odd then $ 20 $ is a quadratic residue modulo $ 23 $. But, neither $ -1 $ nor $ 20 $ are quadratic residues modulo $ 23 $. Thus, $ e^{\Gamma_{5}}\neq 1 $.

Then, we apply Theorem \ref{Matveev11} on the left-hand side of the inequalities \eqref{kakawu2} with the data
\begin{eqnarray*}
&t:=4, \quad \eta_1:=2a, \quad \eta_2:=\alpha,\quad \eta_3:=1+\alpha^{m_1-n_1}, \quad \eta_4:=1+\alpha^{m_2-n_2},\\
&b_1:=k_2-k_1, \quad b_2:=k_2n_1-k_1n_2, \quad b_3:=k_2, \quad b_4:=k_1.
\end{eqnarray*}
Together with combining the right-hand side of \eqref{kakawu2} with the inequalities \eqref{kabaya123} and \eqref{case2}, Theorem \ref{Matveev11} gives
\begin{eqnarray}
n_1&<&3.02\times 10^{16}(n_1-m_1)(n_2-m_2)(\log n_2)\nonumber\\
&<&8.33\times 10^{52}(\log n_2)^{4}.\label{kalo}
\end{eqnarray}
In the above we used the facts that
\begin{eqnarray*}
\min_{1\le i\le 2}\{n_i-m_i\}<1.56\times 10^{11}\log n_2 \quad \text{and}\quad \max_{1\le i\le 2}\{n_i-m_i\}<1.78\times 10^{25}(\log n_2)^{2}.
\end{eqnarray*}
This was obtained under the assumption that the inequality \eqref{kakawu3} does not hold. If \eqref{kakawu3} holds, then so does \eqref{kalo}. Thus, we have that inequality \eqref{kalo} holds provided that inequality \eqref{case2} holds. Otherwise, inequality \eqref{case1} holds which is a better bound than \eqref{kalo}. Hence, conclude that \eqref{kalo} holds in all posibble cases.

By the inequality \eqref{kabizi},
\begin{eqnarray*}
\log\delta \le k_1\log\delta \le n_1\log\alpha +\log 6 <2.38\times 10^{52}(\log n_2)^{4}.
\end{eqnarray*}
By substituting this into \eqref{lemmata1} we get $n_2<4.64\times 10^{137}(\log n_2)^{10}$, and then, by Lemma \ref{gl}, with the data $ r:=10, ~H:=4.64\times 10^{137}$  and $L:=n_2 $, we get that  $ n_2< 4.87\times 10^{165}$. This immediately gives that $ n_1<1.76\times 10^{63} $.

We record what we have proved.
\begin{lemma}\label{lemmata2}
Let $(k_i, n_i, m_i)$ be a solution to $ x_{k_i}=P_{n_i}+P_{m_i} $, with $ 3\le m_i<n_i $ for $ i\in\{1,2\} $ and $ 1\le k_1<k_2 $, then
\begin{eqnarray*}
\max\{k_1, m_1\}<n_1<1.76\times 10^{63}, \quad \text{and} \quad \max\{k_2, m_2\}<n_2<4.87\times 10^{165}.
\end{eqnarray*}
\end{lemma}
\section{Reducing the bounds for $ n_1 $ and $ n_2 $ (I)}
In this section we reduce the bounds for $ n_1 $ and $ n_2 $ given in Lemma \ref{lemmata2} to cases that can be computationally treated. For this, we return to the inequalities for $ \Gamma_3 $, $ \Gamma_4 $ and $ \Gamma_5 $. 
\subsection{The first reduction (I)}
We divide through both sides of the inequality \eqref{basaja1} by $(k_2-k_1)\log\alpha$. We get that
\begin{eqnarray}\label{kapa1}
\left|\dfrac{\log (2a)}{\log\alpha}-\dfrac{k_2n_1-k_1n_2}{k_2-k_1}\right|<\dfrac{36n_2}{\alpha^{\lambda}(k_2-k_1)} \quad \text{with} \quad \lambda: = \min_{1 \le i \le 2}\{n_i-m_i\}.
\end{eqnarray}
We assume that $ \lambda\ge 10 $. Below we apply Lemma \ref{Legendre}. We put $ \tau:= \frac{\log (2a)}{\log\alpha} $, which is irrational and compute its continued fraction \[ [a_0, a_1, a_2, \ldots]=[1, 3, 3, 1, 11, 1, 2, 1, 1, 1, 3, 1, 1, 1, 2, 5, 1, 15, 2, 19, 1, 1, 2, 2, \ldots] \] and its convergents \[ \left[\frac{p_0}{q_0},\frac{p_1}{q_1}, \frac{p_2}{q_2}, \ldots\right] =\left[ 1, \frac{4}{3}, \frac{13}{10}, \frac{17}{13}, \frac{200}{153}, \frac{217}{166}, \frac{634}{485}, \frac{851}{651}, \frac{1485}{1136},
\frac{2336}{1787}, \frac{8493}{6497},  \ldots \right]. \]
Furthermore, we note that taking $ M:=4.87\times 10^{165} $ (by Lemma \ref{lemmata2}), it follows that
\begin{eqnarray*}
q_{315}>M>n_2>k_2-k_1 \quad \text{and}\quad a(M):=\max\{a_i:0\le i\le 315\}=a_{282}=2107.
\end{eqnarray*}
Thus, by Lemma \ref{Legendre}, we have that
\begin{eqnarray}\label{kapa2}
\left|\tau - \dfrac{k_2n_1-k_1n_2}{k_2-k_1}\right|>\dfrac{1}{2109(k_2-k_1)^{2}}.
\end{eqnarray}
Hence, combining the inequalities \eqref{kapa1} and \eqref{kapa2}, we obtain
\begin{eqnarray*}
\alpha^{\lambda}<75924n_2(k_2-k_1)<1.75\times 10^{336},
\end{eqnarray*}
so $ \lambda\le 2714 $. This was obtained under the assumption that $ \lambda\ge 10 $, Otherwise, $ \lambda<10<2714 $ holds as well.

Now, for each $ n_i-m_i = \lambda\in [1, 2714] $ we estimate a lower bound $ |\Gamma_4| $, with
\begin{eqnarray}\label{LLL1}
\Gamma_4&=&(k_j-k_i)\log(2a)+(k_jn_i-k_in_j)\log\alpha+k_j\log(1+\alpha^{m_i-n_i})
\end{eqnarray}
given in the inequality \ref{kipro1}, via the procedure described in Subsection \ref{Reduction}  (LLL-algorithm). We recall that $ \Gamma_4\neq 0 $. 

We apply Lemma \ref{LLL} with the data: 
\begin{eqnarray*}
&t:=3, \quad \tau_1:=\log(2a), \quad \tau_2:=\log\alpha, \quad \tau_3:=\log(1+\alpha^{-\lambda}),\\
&x_1:=k_j-k_i, \quad x_2:=k_jn_i-k_in_j, \quad x_3:=k_j.
\end{eqnarray*}
We set $ X:= 5.4\times 10^{166} $ as an upper bound to $ |x_i|<11n_2 $ for all $ i=1, 2, 3 $, and $ C:=(20X)^{5} $. A computer in \textit{Mathematica} search allows us to conclude, together with the inequality \eqref{kipro1}, that
\begin{eqnarray*}
2\times 10^{-671}<\min_{1\le \lambda \le  2714}|\Gamma_4|<8n_2\alpha^{-\nu},\quad \text{with} \quad \nu:=\min\{n_i, n_j-m_j\}
\end{eqnarray*}
which leads to $ \nu \le 6760 $. As we have noted before, $ \nu = n_1 $ (so $n_1\le 6760$) or $ \nu = n_j-m_j $.

Next, we suppose that $ n_j-m_j = \nu \le 6760 $. Since $ \lambda\le 2714 $, we have
\begin{eqnarray*}
\lambda := \min_{1\le i \le 2}\{n_i-m_i\}\le 2714 \quad \text{and} \quad \chi :=\max_{1\le i \le 2}\{n_i-m_i\}\le 6760.
\end{eqnarray*}
Now, returning to the inequality \eqref{kakawu1} which involves
\begin{eqnarray}\label{LLL2}
\Gamma_5:&=&(k_2-k_1)\log(2a)+(k_2n_1-k_1n_2)\log\alpha \nonumber\\&&+k_2\log (1+\alpha^{m_1-n_1})-k_1\log (1+\alpha^{m_2-n_2})\neq 0,
\end{eqnarray}
we use again the LLL-algorithm to estimate the lower bound for $ |\Gamma_5| $ and thus, find a bound for $ n_1 $ that is better than the one given in Lemma \ref{lemmata2}.

We distinguish the cases $ \lambda < \chi $ and $ \lambda = \chi $. 
\subsection{The case $ \lambda < \chi $}
We take $ \lambda \in [1, 2714] $ and $ \chi \in [\lambda+1, 6760] $ and apply Lemma \ref{LLL} with the data:
\begin{eqnarray*}
& t:=4,\quad \tau_1:=\log(2a), \quad \tau_2:= \log\alpha, \quad \tau_3: = \log(1+\alpha^{m_1-n_1}), \quad  \tau_4: = \log(1+\alpha^{m_2-n_2}),\\
& x_1:=k_2-k_1, \quad x_2:= k_2n_1-k_1n_2, \quad x_3: = k_2, \quad x_4:=-k_1.
\end{eqnarray*}
We also put $ X:=5.4\times 10^{166} $ and $ C:=(20X)^{9} $. After a computer search in \textit{Mathematica} together with the inequality \ref{kakawu1}, we can confirm that
\begin{eqnarray}
8\times 10^{-1342}<\min_{\substack{1\le \lambda\le 2714 \\ \lambda+1\le \chi \le 6760}}|\Gamma_5| < 6n_2 \alpha^{-n_1}.
\end{eqnarray}
This leads to the inequality
\begin{eqnarray}
\alpha^{n_1}< 7.5\times 10^{1341}n_2.
\end{eqnarray}
Subsitituting for the bound  $ n_2 $ given in Lemma \ref{lemmata2}, we get that $ n_1\le 12172 $.
\subsection{The case $ \lambda=\chi $}
In this case, we have 
\begin{eqnarray*}
\Lambda_5:=(k_2-k_1)(\log(2a)+\log(1+\alpha^{m_1-n_1}))+(k_2n_1-k_1n_2)\log\alpha \neq 0.
\end{eqnarray*}
We divide through the inequality \ref{kakawu1} by $(k_2-k_1)\log\alpha$ to obtain
\begin{eqnarray}\label{kaka1}
\left|\dfrac{\log(2a)+\log(1+\alpha^{m_1-n_1})}{\log \alpha}-\dfrac{k_2n_1-k_1n_2}{k_2-k_1}\right|<\dfrac{21n_2}{\alpha^{n_1}(k_2-k_1)}
\end{eqnarray}
We now put $$ \tau_{\lambda}:=\dfrac{\log(2a)+\log(1+\alpha^{-\lambda})}{\log \alpha} $$ and compute its continued fractions $ [a_0^{(\lambda)}, a_1^{(\lambda)} , \ldots] $ and its convergents $[p_0^{(\lambda)}/q_0^{(\lambda)}, p_1^{(\lambda)}/q_1^{(\lambda)}, \ldots]$ for each $ \lambda \in [1, 2714] $. Furthermore, for each case we find an integer $ t_{\lambda} $ such that $ q_{t_{\lambda}}^{(\lambda)}>M:=4.87\times 10^{165}>n_2>k_2-k_1 $ and calculate
\begin{eqnarray*}
a(M):=\max_{1\le\lambda\le 2714}\left\{a_{i}^{(\lambda)}: 0 \le i \le t_{\lambda}\right\}.
\end{eqnarray*}
A computer search in \textit{Mathematica} reveals that for  $ \lambda = 321 $, $ t_{\lambda} = 330  $ and $ i=263 $, we have that $ a(M) = a_{321}^{(330)}=306269 $. Hence, combining the conclusion of Lemma \ref{Legendre} and the inequality \eqref{kaka1}, we get
\begin{eqnarray*}
\alpha^{n_1}< 21\times 306271n_2(k_2-k_1)< 1.525\times 10^{338},
\end{eqnarray*}
so $ n_1\le 2730 $. Hence, we obtain that $ n_1\le 12172 $ holds in all cases ($\nu =n_1$, $\lambda < \chi$ or $\lambda = \chi$). By the inequality \eqref{kabizi}, we have that
\begin{eqnarray*}
\log\delta\le k_1\log\delta \le n_1\log\alpha +\log 6 <3475.
\end{eqnarray*}
By considering the second inequality in \eqref{lemmata1}, we can conclude that $ n_2\le 9.9\times 10^{39}(\log n_2)^2 $, which immediately yields $ n_2<3.36\times 10^{44}  $, by a simple application of Lemma \ref{gl}. We summarise the first cycle of our reduction process as follows:
\begin{eqnarray}
n_1\le 12172 \quad \text{and} \quad n_2\le 3.36\times 10^{44}.
\end{eqnarray}
From the above, we note that the upper bound on $ n_2 $ represents a very good reduction of the bound given in Lemma \ref{lemmata2}. Hence, we expect that if we restart our reduction cycle with the new bound on $ n_2 $, then we get a better bound on $ n_1 $. Thus, we return to the inequality \eqref{kapa1} and take $ M:=3.36\times 10^{44} $. A computer search in \textit{Mathematica} reveals that
\begin{eqnarray*}
q_{88}>M>n_2>k_2-k_1 \quad \text{and} \quad a(M):=\max\{a_i: 0\le i\le 88\}=a_{54} =373,
\end{eqnarray*}
from which it follows that $ \lambda \le 752 $. We now return to \eqref{LLL1} and we put $ X:=3.36\times 10^{44} $ and $ C:=(10X)^{5} $ and then apply the LLL-algorithm in Lemma \ref{LLL} to $ \lambda \in[1, 752] $. After a computer search, we get
\begin{eqnarray*}
5.33\times 10^{-184}<\min_{1\le \lambda\le 752}|\Gamma_4| < 8n_2\alpha^{-\nu},
\end{eqnarray*}
then $ \nu \le 1846 $. By continuing under the assumption that $ n_j-m_j=\nu \le 1846 $, we return to \eqref{LLL2} and put $ X:=3.36\times 10^{44} $, $ C:=(10X)^9 $ and $ M:=3.36\times 10^{44} $ for the case $ \lambda <\chi $ and $ \lambda = \chi $. After a computer search, we confirm that
\begin{eqnarray*}
2\times 10^{-366}< \min_{\substack{1\le \lambda \le 752\\ \lambda+1\le \chi \le 1846}}|\Gamma_5|<6n_2\alpha^{-n_1},
\end{eqnarray*}
gives $ n_1\le 3318 $, and
$
a(M)=a_{175}^{(205)}=206961
$, leads to $ n_1\le 772 $. Hence, in both cases $ n_1\le 3318 $ holds. This gives $ n_2\le 5\times 10^{42} $ by a similar procedure as before, and $ k_1\le  $.

We record what we have proved.

\begin{lemma}\label{firstredn1}
Let $ (k_i, n_i, m_i) $ be a solution to $ X_i=P_{n_i}+P_{m_i} $, with $ 3\le m_i<n_i $ for $ i=1,2 $ and $ 1\le k_1<k_2 $, then
\begin{eqnarray*}
m_1<n_1\le 3318, \quad  k_1\le 3125 \quad \text{and} \quad n_2\le 5\times 10^{42}.
\end{eqnarray*}
\end{lemma}
\subsection{The final reduction (I)}
Returning back to \eqref{Kay1} and \eqref{Kayii1} and using the fact that $ (x_1, y_1) $ is the smallest positive solution to the Pell equation \eqref{Pelleqn1}, we obtain
\begin{eqnarray*}
x_k &=& \dfrac{1}{2}(\delta^k+\sigma^k) \quad =\quad \dfrac{1}{2}\left(\left(x_1+y_1\sqrt{d}\right)^{k}+\left(x_1-y_1\sqrt{d}\right)^k\right)\\
&=&\dfrac{1}{2}\left(\left(x_1+\sqrt{x_1^2\mp 1}\right)^{k}+\left(x_1-\sqrt{x_1^2\mp 1}\right)^k\right): = Q^{\pm}_{k}(x_1).
\end{eqnarray*}
Thus, we return to the Diophantine equation $ x_{k_1}=P_{n_1}+P_{m_1} $ and consider the equations 
\begin{eqnarray}\label{water1}
Q^{+}_{k_1}(x_1)=P_{n_1}+P_{m_1} \quad \text{and} \quad Q^{-}_{k_1}(x_1)=P_{n_1}+P_{m_1},
\end{eqnarray}
with $ k_1\in [1, 3125] $, $ m_1\in [3,3318 ] $ and $n_1\in[m_1+1, 3318 ] $.

Besides the trivial case $ k_1=1 $, with the help of a computer search in \textit{Mathematica} on the above equations in \eqref{water1}, we list the only nontrivial solutions in the tables below. We also note that $ 3+2\sqrt{2}=(1+\sqrt{2})^2 $, so these solutions come from the same Pell equation when $ d=2 $.

\begin{table}[H]
\parbox{.40\linewidth}{
\centering
\begin{tabular}{|c c c c c|}
\hline
& & &$ Q^{+}_{k_1}(x_1) $&\\
\hline
$ k_1 $& $ x_1 $& $ y_1 $& $ d $ & $ \delta $\\
\hline
$ 2 $&$ 2 $& $ 1 $&$ 3 $ & $ 2+\sqrt{3} $\\
$ 2 $&$ 3 $& $ 2 $&$ 2 $ & $ 3+2\sqrt{2} $\\
$ 2 $&$ 4 $& $ 1 $&$ 15 $ & $ 4+\sqrt{15} $\\
$ 2 $&$ 5 $& $ 2 $&$ 6 $ & $ 5+2\sqrt{6} $\\
$ 2 $&$ 21 $& $ 2 $&$ 110 $ & $ 21+2\sqrt{110} $\\
$ 2 $&$ 22 $& $ 1 $&$ 483 $ & $ 22+\sqrt{483} $\\
$ 2 $&$ 47 $& $ 4 $&$ 138 $ & $ 47+4\sqrt{138} $\\
\hline
\end{tabular}
}
\hfill
\parbox{.40\linewidth}{
\centering
\begin{tabular}{|c c c c c|}
\hline
& & &$ Q^{-}_{k_1}(x_1) $&\\
\hline
$ k_1 $& $ x_1 $& $ y_1 $& $ d $ & $ \delta $\\
\hline
$ 2 $&$ 1 $& $ 1 $&$ 2 $ & $ 1+\sqrt{2} $\\
$ 2 $&$ 2 $& $ 1 $&$ 5 $ & $ 2+\sqrt{5} $\\
$ 2 $&$ 3 $& $ 1 $&$ 10 $ & $ 3+\sqrt{10} $\\
$ 2 $&$ 4 $& $ 1 $&$ 17 $ & $ 4+\sqrt{17} $\\
$ 2 $&$ 5 $& $ 1 $&$ 26 $ & $ 5+\sqrt{26} $\\
$ 2 $&$ 9 $& $ 1 $&$ 82 $ & $ 9+\sqrt{82} $\\
$ 2 $&$ 10 $& $ 1 $&$ 101 $ & $ 10+\sqrt{101} $\\
$ 2 $&$ 17 $& $ 1 $&$ 290 $ & $ 17+\sqrt{290} $\\
$ 2 $&$ 42 $& $ 1 $&$ 1765 $ & $ 42+\sqrt{1765} $\\
$ 2 $&$ 47 $& $ 1 $&$ 2210 $ & $ 47+\sqrt{2210} $\\
$ 2 $&$ 63 $& $ 1 $&$ 3970 $ & $ 63+\sqrt{3970} $\\
\hline
\end{tabular}
}
\end{table}
From the above tables, we set each $ \delta:=\delta_{t} $ for $ t=1, 2, \ldots 17 $.
We then work on the linear forms in logarithms $ \Gamma_1 $ and $ \Gamma_2 $, in order to reduce the bound on $ n_2 $ given in  Lemma \ref{firstredn1}. From the inequality \eqref{mukazi1}, for $ (k,n,m):=(k_2, n_2, m_2) $, we write
\begin{eqnarray}\label{mukazi111}
\left|k_2\dfrac{\log\delta_t}{\log\alpha}-n_2+\dfrac{\log(2a)}{\log (\alpha^{-1})}\right|<\left(\frac{5}{\log \alpha}\right)\alpha^{-(n_2-m_2)}, \quad \text{for}\quad t=1,2, \ldots 17.
\end{eqnarray}
We put
\begin{eqnarray*}
\tau_{t}:=\dfrac{\log\delta_t}{\log\alpha}, \qquad \mu_t:=\dfrac{\log(2a)}{\log(\alpha^{-1})}\qquad \text{and} \quad (A_t, B_t):=\left(\frac{5}{\log \alpha}, \alpha\right).
\end{eqnarray*}
We note that $ \tau_t $ is transcendental by the Gelfond-Schneider's Theorem and thus, $ \tau_t $ is irrational. We can rewrite the above inequality, \ref{mukazi111} as
\begin{eqnarray}\label{mukaziii11}
0<|k_2\tau_t-n_2+\mu_t|<A_tB_t^{-(n_2-m_2)}, \quad \text{for} \quad t=1, 2, \ldots, 17.
\end{eqnarray}
We take $ M:= 5\times 10^{42} $ which is the upper bound on $ n_2 $ according to Lemma \ref{firstredn1} and apply Lemma \ref{Dujjella} to the inequality \eqref{mukaziii11}. As before, for each $ \tau_t $ with $ t=1, 2, \ldots, 17 $, we compute its continued fraction $ [a_0^{(t)}, a_1^{(t)}, a_2^{(t)}, \ldots ] $ and its convergents $ p_0^{(t)}/q_0^{(t)}, p_1^{(t)}/q_1^{(t)}, p_2^{(t)}/q_2^{(t)}, \ldots $. For each case, by means of a computer search in \textit{Mathematica}, we find and integer $s_{t}$ such that
\begin{eqnarray*}
q^{(t)}_{s_t}> 3\times 10^{43}=6M \qquad \text{ and } \qquad \epsilon_t:=||\mu_t q^{(t)}||-M||\tau_t q^{(t)}|>0.
\end{eqnarray*}
We finally compute all the values of $ b_t:=\lfloor \log(A_t q^{(t)}_{s_t}/\epsilon_t)/\log B_t \rfloor $. The values of  $ b_t $ correspond to the upper bounds on $ n_2-m_2 $, for each $ t=1, 2, \ldots, 17 $, according to Lemma \ref{Dujjella}. The results of the computation for each $ t $ are recorded in the table below.

\vspace*{0.5cm}
\begin{center}
\begin{tabular}{lllccc}
\hline 
$t$ & $\delta_t$ & $s_t$ & $q_{s_t} $ & $\epsilon_t>$ & $b_t$ \\ 
\hline 
$1$ & $2+\sqrt{3}$ & $85$ & $ 8.93366\times 10^{43} $ & $ 0.3100 $ & $374$ \\ 
$2$ & $ 4+\sqrt{15} $ & 90 & $ 3.90052\times 10^{43} $ & $ 0.3124 $ & $371$ \\ 
$3$ & $ 5+2\sqrt{6} $ & $ 80 $ & $ 3.16032\times 10^{43} $ & $ 0.0122 $ & $382$ \\  
$4$ & $ 21+2\sqrt{110} $& $ 88 $ & $ 6.33080\times 10^{43} $ & $ 0.2200 $ & $374$ \\ 
$5$ & $22+\sqrt{483}$ & $ 75 $ & $4.19689\times 10^{43}$ & $0.2361$ & $372$ \\  
$6$ & $47+4\sqrt{138} $ & $ 96 $ & $ 7.76442\times 10^{43} $ & $ 0.3732 $ & $373$ \\  
$7$ & $1+\sqrt{2} $ & $ 78 $ & $  1.46195 \times 10^{44}$ & $ 0.3328 $ & $375$ \\  
$8$ & $2+\sqrt{5} $ & $ 94 $ & $ 1.48837\times 10^{44} $ & $ 0.2146 $ & $377$ \\ 
$9$ & $3+\sqrt{10}$  & $ 88 $ & $ 4.21425\times 10^{43} $ & $ 0.1347 $& $374$ \\  
$10$ & $4+\sqrt{17}$ & $ 92 $ & $ 1.11753\times 10^{44} $ & $ 0.2529 $ & $375$ \\ 
$11$ & $5+\sqrt{26}$ & $ 98 $ & $ 3.23107\times 10^{43} $ & $ 0.1043 $ & $374$ \\  
$12$ & $ 9+\sqrt{82} $ & $ 74 $ & $ 5.25207\times 10^{43} $ & $ 0.2181 $ & $373$ \\  
$13$ & $10+\sqrt{101} $ & $ 94 $ & $ 1.86122\times 10^{44} $ & $ 0.2672 $ & $377$ \\  
$14$ & $17+\sqrt{290} $ & $ 87 $ & $ 1.06422\times 10^{44} $ & $ 0.0193 $ & $384$ \\  
$15$ & $42+\sqrt{1765} $ & $ 78 $ & $ 3.81406\times 10^{43} $ & $ 0.1768 $ & $373$ \\  
$16$ & $47+\sqrt{2210} $ & $ 94 $ & $ 3.92482\times 10^{43} $ & $ 0.4476 $ & $370$ \\ 
$17$ & $63+\sqrt{3970}$ & $85$ & $6.00550\times 10^{43}$ & $ 0.4056 $ & $371$ \\ 
\hline 
\end{tabular} 
\end{center}
\vspace*{0.5cm}

By replacing $ (k, n, m):=(k_2, n_2, m_2) $ in the inequality \eqref{kabuzi}, we can write
\begin{eqnarray}\label{water112}
\left|k_2\dfrac{\log\delta_t}{\log\alpha}-n_2+\dfrac{\log(2a(1+\alpha^{-(n_2-m_2)}))}{\log(\alpha^{-1})}\right|< \left(\dfrac{3}{\log\alpha}\right)\alpha^{-n_2}, 
\end{eqnarray}
for  $t=1,2,\ldots, 17$.

We now put
\begin{eqnarray*}
\tau_{t}:=\dfrac{\log\delta_t}{\log\alpha}, \quad \mu_{t, n_2-m_2}:=\dfrac{\log(2a(1+\alpha^{-(n_2-m_2)}))}{\log(\alpha^{-1})}\quad \text{and} \quad (A_t, B_t):=\left(\frac{3}{\log \alpha}, \alpha\right).
\end{eqnarray*}
With the above notations, we can rewrite \eqref{water112} as
\begin{eqnarray}\label{water21}
0<|k_2\tau_t - n_2+\mu_{t, n_2-m_2}|<A_tB_t^{-n_2}, \quad \text{ for} \quad t=1,2, \ldots 17.
\end{eqnarray}
We again apply Lemma \ref{Dujjella} to the above inequality \eqref{water21}, for 
\begin{eqnarray*}
t=1, 2, \ldots, 17, \quad n_2-m_2 =1, 2, \ldots, b_t, \quad \text{with}\quad M:=5\times 10^{43}.
\end{eqnarray*}
We take
\begin{eqnarray*}
\epsilon_{t, n_2-m_2}:=||\mu_t q^{(t, n_2-m_2)}|| -M||\tau_t q^{(t, n_2-m_2)}||>0,
\end{eqnarray*}
and 
\begin{eqnarray*}
b_t=b_{t, n_2-m_2}:=\lfloor \log(A_t q^{(t, n_2-m_2)}_{s_t}/\epsilon_{t, n_2-m_2})/\log B_t \rfloor.
\end{eqnarray*}
With the help of \text{Mathematica}, we obtain that
\vspace{0.4cm}

\begin{tabular}{c|ccccccccccccc}
\hline
$t$& $ 1 $&$ 2 $&$ 3 $&$ 4 $&$ 5 $&$ 6 $&$ 7 $&$ 8 $&$ 9 $\\

$b_{t, n_2-m_2}$& $ 388 $&$ 389 $&$ 394 $&$ 394 $&$ 393 $&$ 394 $&$ 396 $&$ 392 $&$ 392$\\
\hline
$t$&$ 10 $&$ 11 $&$ 12 $&$ 13 $&$ 14 $ &$ 15 $&$ 16 $&$ 17 $\\

$b_{t, n_2-m_2}$& $396 $&$ 392 $&$ 408 $&$ 390 $&$ 396 $&$ 396 $& $ 388 $&$ 389 $\\
\hline
\end{tabular}

\begin{eqnarray*}
\text{Thus, }~\max\{b_{t, n_2-m_2}: t=1, 2, \ldots, 17 \quad \text{and} \quad n_2-m_2 = 1, 2, \ldots b_t\} \le 408 .
\end{eqnarray*}
Thus, by Lemma \ref{Dujjella}, we have that $ n_2\le 408 $, for all $ t=1,2, \ldots, 17 $, and by the inequality \eqref{faji1} we have that $ n_1\le n_2+4 $. From the fact that $ \delta^{k}\le 2\alpha^{n+3} $, we can conclude that $ k_1< k_2\le 133  $. Collecting everything together, our problem is reduced to search for the solutions for \eqref{kalayim2} in the following range
\begin{eqnarray}
1\le k_1<k_2\le 133, \quad 0\le m_1<n_1 \in [3, 408] \quad \text{and} \quad 0\le m_2<n_2 \in [3, 408].
\end{eqnarray}
After a computer search on the equation \eqref{kalayim2} on the above ranges, we obtained the following solutions, which are the only solutions for the  exceptional $ d $ cases  we have stated in Theorem \ref{Main1}:

For the  $ +1 $ case:
\begin{eqnarray*}
(d=2)&&x_1 =3=P_6+P_0=P_5+P_3, \quad x_2=17=P_{12}+P_3;\\
(d=3)&&x_1 =2=P_3+P_0=P_3+P_3,~   x_2=7=P_9+P_0=P_7+P_6,~ \\&&  x_3=26=P_{13}+P_8;\\
  (d=6)&&x_1 =5=P_8+P_0 = P_7+P_3= P_6+P_5,\\&& x_2=49=P_{16}+P_0=P_{15}+P_{12}=P_{14}+P_{13}; \\
  (d=15)&&x_1=4=P_7+P_0=P_6+P_3=P_5+P_5, \quad x_2=31=P_14+P_6;\\
 (d=110)&& x_1=21=P_{13}+P_0=P_{12}+P_{8}=P_{11}+P_{10}, \\ && x_2=881=P_{26}+P_{17}=P_{25}+P_{22};\\
 (d=483)&& x_1=22=P_{13}+P_{3}, \quad x_2=967=P_{26}+P_{20}=P_{25}+P_{23}.
\end{eqnarray*}

For the $ -1 $ case:
\begin{eqnarray*}
(d=2) && x_1=1=P_3+P_0, \quad x_2=7=P_9+P_0=P_8+P_5=P_7+P_6,\\
&& x_3=41=P_{15}+P_7=P_{14}+P_{10}=P_{13}+P_{12};\\
(d=5)&& x_1=2=P_5+P_0=P_3+P_3, \quad x_2=38=P_{15}+P_{3};\\
(d=10)&& x_1=3=P_6+P_0 = P_5+P_3, \quad x_2=117=P_{19}+P_{6};\\
(d=17)&& x_1=4=P_7+P_0=P_6+P_3=P_5+P_5,\quad x_2=P_{22}+P_{6}.
\end{eqnarray*}
This completes the proof of Theorem \ref{Main1}. \qed

\section{Proof of Theorem \ref{Main2}}
The proof of Theorem \ref{Main2} will be similar to that of Theorem \ref{Main1}. We also give the details for the benefit of the reader. Further, for technical reasons in our proof, we assume that $ d\ge 5 $ and then treat the cases $ d\in \{2,3\} $ during the reduction procedure.

Let $ (X_1, Y_1) $ be the smallest positive integer solution to the Pell quation \eqref{Pelleqn2}. We Put
\begin{eqnarray}\label{kay2}
\rho:=\dfrac{X_1+Y_1\sqrt{d}}{2} \quad \text{and} \quad \varrho=\dfrac{X_1-Y_1\sqrt{d}}{2}.
\end{eqnarray}
From which we get that
\begin{eqnarray}\label{kalayi2}
\rho\cdot\varrho=\dfrac{X_1^2-dY_1^2}{4} =: \epsilon, \quad \text{where} \quad \epsilon\in\{\pm 1\}. 
\end{eqnarray}
Then
\begin{eqnarray}\label{kayiii2}
X_n = \rho^k+\varrho^k.
\end{eqnarray}
Since $ \rho \ge \frac{1+\sqrt{5}}{2} $, it follows that the estimate 
\begin{eqnarray}\label{kalayim11}
\dfrac{\rho^{k}}{\alpha^{2}}\le X_k \le 2\rho^k \quad \text{ holds for all } \quad k\ge 1.
\end{eqnarray}
Similarly, as before, we assume that $ (k_1, n_1, m_1) $ and $ (k_2, n_2, m_2) $ are triples of integers such that
\begin{eqnarray}\label{kalayim12}
X_{k_1}=P_{n_1}+P_{m_1} \quad \text{and} \quad X_{k_2}=P_{n_2}+P_{m_2}
\end{eqnarray}
We asuume that $ 1\le k_1 < k_2 $. We also assume that $ 4\le m_j< n_j $ for $ j=1,2 $.
We set $ (k,n,m):=(k_j, n_j, m_j) $, for $ j=1,2 $. Using the inequalities \eqref{Pado6} and \eqref{kalayim11}, we get from \eqref{kalayim12} that 
\begin{eqnarray*}
\dfrac{\rho^{k}}{\alpha^{2}}\le X_k=P_n+P_m\le 2\alpha^{n-1}\quad \text{and} \quad \alpha^{n-2}\le P_n+P_m = X_k\le 2\rho^{k}.
\end{eqnarray*}
The above inequalities give
\begin{eqnarray*}
(n-2)\log\alpha-\log2 <k\log\rho < (n+1)\log\alpha +\log 2.
\end{eqnarray*}
Dividing through by $ \log\alpha $ and setting $ c_1:=1/\log\alpha $, as before, we get that 
\begin{eqnarray*}
-2-c_1\log 2 <c_1k\log\rho-n<1+c_1\log2,
\end{eqnarray*}
and since $ \alpha^3>2 $, we get
\begin{eqnarray}\label{kabizi2}
|n-c_1\log\rho|<5.
\end{eqnarray}
Furthermore, $ k<n $, for if not, we would then get that 
\begin{eqnarray*}
\rho^{n}\le \rho^{k}<2\alpha^{n+1}, \quad \text{implying}\quad \left(\dfrac{\rho}{\alpha}\right)^{n}<2\alpha,
\end{eqnarray*}
which is false since $ \rho\le \frac{1+\sqrt{5}}{2} $, $ 1.32<\alpha<1.33 $ and $ n\ge 5 $.

Besides, given that $ k_1<k_2 $, we have by \eqref{Pado7} and \eqref{kalayim12} that
\begin{eqnarray*}
\alpha^{n_1-2}\le P_{n_1}\le P_{n_1}+P_{m_1}=X_{k_1}<X_{k_2}= P_{n_2}+P_{m_2} \le 2P_{n_2}<2\alpha^{n_2-1}.
\end{eqnarray*}
Thus, as before, we get that
\begin{eqnarray}\label{faji2}
n_1<n_2+4.
\end{eqnarray}
\subsection{An inequality for $ n $ and $ k $ (II)}
Using the equations \eqref{Pado3} and \eqref{kay2} and \eqref{kalayim12}, we get
\begin{eqnarray*}
\rho^{k}+\varrho^{k}=P_n+P_m=a\alpha^{n}+e(n)+a\alpha^{m}+e(m)
\end{eqnarray*}
So,
\begin{eqnarray*}
\rho^{k}-a(\alpha^{n}+\alpha^{m})=-\varrho^{k}+e(n)+e(m),
\end{eqnarray*}
and by \eqref{Pado6}, we have
\begin{eqnarray*}
\left|\rho^{k}a^{-1}\alpha^{-n}(1+\alpha^{m-n})^{-1}-1\right|&\leq& \dfrac{1}{\rho^{k}a(\alpha^{n}+\alpha^{m})}+\dfrac{2|b|}{\alpha^{n/2}a(\alpha^{n}+\alpha^{m})}\\&&+\dfrac{2|b|}{\alpha^{m/2}a(\alpha^{n}+\alpha^{m})}\\
&\le& \dfrac{1}{a\alpha^{n}}\left(\dfrac{1}{\rho^{k}}+\dfrac{2|b|}{\alpha^{n/2}}+\dfrac{2|b|}{\alpha^{m/2}}\right)<\dfrac{2.5}{\alpha^{n}}.
\end{eqnarray*}
Thus, we have
\begin{eqnarray}\label{kapuli2}
\left|\rho^{k}a^{-1}\alpha^{-n}(1+\alpha^{m-n})^{-1}-1\right|&<&\dfrac{2.5}{\alpha^{n}}.
\end{eqnarray}
Put
\begin{eqnarray*}
\Lambda_1^{'}:=\rho^{k}a^{-1}\alpha^{-n}(1+\alpha^{m-n})^{-1}-1
\end{eqnarray*}
and 
\begin{eqnarray*}
\Gamma_1^{\prime}:=k\log\rho-\log a -n\log\alpha-\log(1+\alpha^{m-n}).
\end{eqnarray*}
Since $|\Lambda_1^{\prime}|=|e^{\Gamma_1^{\prime}}-1|<0.83$ for $ n\ge 4 $ (because $2.5/\alpha^{4}<0.83$), it follows that $ e^{|\Gamma_1^{\prime}|}<4 $ and so
\begin{eqnarray*}
|\Gamma_1^{\prime}|<e^{|\Gamma_1^{\prime}|}|e^{\Gamma_1^{\prime}}-1|<\dfrac{10}{\alpha^{n}}.
\end{eqnarray*}
Thus, we get that
\begin{eqnarray}\label{kabuzi22}
\left|k\log\rho-\log a -n\log\alpha-\log(1+\alpha^{m-n})\right|<\dfrac{10}{\alpha^{n}}.
\end{eqnarray}
We apply Theorem \ref{Matveev11} on the left-hand side of \eqref{kapuli2} with the data:
\begin{eqnarray*}
&t:=4, \quad \eta_{1}:=\rho, \quad \eta_2:=a, \quad \eta_3:=\alpha, \quad \eta_4: =1+\alpha^{m-n},\\
&b_1:=k, \quad b_2:=-1, \quad b_3:=-n, \quad b_4:=-1.
\end{eqnarray*}
Furthermore, we take same the number field as before, $ \mathbb{K}=\mathbb{Q}(\sqrt{d}, \alpha) $ with degree $ D=6 $. We also  take $D_{\mathbb{K}}=n$.
First we note that the left-hand side of \eqref{kapuli1} is non-zero, since otherwise,
\begin{eqnarray*}
\rho^{k}=a(\alpha^{n}+\alpha^{m}).
\end{eqnarray*}
By the same argument as before, we get a contradiction. Thus, $ \Lambda_1^{\prime}\neq 0 $ and we can apply Theorem \ref{Matveev11}. Further, $$ a=\dfrac{\alpha(\alpha+1)}{3\alpha^2-1}, $$ the mimimal polynomial of $ a $ is $ 23x^3-23x^2+6x-1 $ and has roots $ a, b, c $. Since $ \max\{a,b,c\}<1 $ (by \eqref{Pado5}), then $h(\eta_2)=h(a)=\frac{1}{3}\log 23$. Thus, we can take $ A_1:=3\log\rho $, $A_2:=2\log 23$,  $A_3:=2\log\alpha$, and  $A_4:=2(n-m)\log\alpha+6\log 2$.

Now, Theorem \ref{Matveev11} tells us that
\begin{eqnarray*}
\log|\Lambda_1^{\prime}|&> &-1.4\times 30^{7}\times 4^{4.5}\times 6^{2}(1+\log 6)(1+\log n)(3\log\rho)\\
&& \times (2\log 23)(2\log\alpha)(2(n-m)\log\alpha+6\log 2)\\
&>&-2.08\times 10^{17}(n-m)(\log n)(\log\rho).
\end{eqnarray*}
Comparing the above inequality with \eqref{kapuli2}, we get
\begin{eqnarray*}
n\log\alpha - \log 2.5 < 2.08\times 10^{17}(n-m)(\log n)(\log\rho).
\end{eqnarray*}
Hence, we get that
\begin{eqnarray}\label{good12}
n<7.40\times 10^{17}(n-m)(\log n)(\log\rho).
\end{eqnarray}
We now return to the equation $ X_k=P_n+P_m $ and rewrite it as
\begin{eqnarray*}
\rho^{k}-a\alpha^{n} = -\varrho^{k}+e(n)+P_m,
\end{eqnarray*}
we obtain
\begin{eqnarray}\label{kayija2}
\left|\rho^{k}a^{-1}\alpha^{-n}-1\right|\leq \dfrac{1}{a\alpha^{n-m}}\left(\dfrac{1}{\alpha}+\dfrac{1}{\alpha^{m+n/2}}+\dfrac{1}{\rho^{k}\alpha^{m}}\right)<\dfrac{3}{\alpha^{n-m}}.
\end{eqnarray}
Put
\begin{eqnarray*}
\Lambda_2^{\prime}:=\rho^{k}a^{-1}\alpha^{-n}-1, \quad \Gamma_{2}^{\prime}:=k\log\rho-\log a-n\log\alpha.
\end{eqnarray*}
We assume for technical reasons that $ n-m\ge 10 $. So $ |e^{\Lambda_2}-1|<\frac{1}{2} $. It follows that 
\begin{eqnarray}\label{mukazi2}
\left|k\log\rho-\log a-n\log\alpha\right|=|\Gamma_2^{\prime}|<e^{|\Lambda_2^{\prime}|}|e^{\Lambda_2^{\prime}}-1|<\dfrac{6}{\alpha^{n-m}}.
\end{eqnarray}
Furthermore, $ \Lambda_2^{\prime}\neq 0 $ (so $\Gamma_2^{\prime} \neq 0$), since $ \rho^{k}\in\mathbb{Q}(\alpha) $ by the previous argument.

We now apply Theorem \ref{Matveev11} to the left-hand side of \eqref{kayija2} with the data
\begin{eqnarray*}
t:=3, \quad \eta_1:=\rho, \quad\eta_2:=a, \quad \eta_3:=\alpha, \quad b_1:=k, \quad b_2:=-1, \quad b_3:=-n.
\end{eqnarray*}
Thus, we have the same $ A_1, ~A_2, A_3 $ as before. Then, by Theorem \ref{Matveev11}, we conclude that
\begin{eqnarray*}
\log|\Lambda|>-9.50\times 10^{14}(\log\rho)(\log n)(\log\alpha).
\end{eqnarray*}
By comparing with \eqref{kayija2}, we get
\begin{eqnarray}\label{good22}
n-m <9.52\times 10^{14}(\log\rho)(\log n).
\end{eqnarray}
This was obtained under the assumption that $ n-m\ge 10 $, but if $ n-m<10$, then the inequality also holds as well. We replace the bound \eqref{good22} on $ n-m $ in \eqref{good12} and use the fact that $ \rho^{k}\le 2\alpha^{n+1} $, to obtain bounds on $ n $ and $ k $ in terms of $ \log n $ and $\log\rho$. We again record  what we have proved.
\begin{lemma}
Let $ (k,n,m) $ be a solution to the equation  $X_k=P_n+P_m$ with $ 3\le m<n $, then
\begin{eqnarray}\label{lemmata1234}
k< 1.98\times 10^{32}(\log n)^{2}(\log\rho) \quad \text{and} \quad n<7.03\times 10^{32}(\log n)^{2}(\log\rho)^{2}.
\end{eqnarray}
\end{lemma}
\subsection{Absolute bounds (II)}
We recall that $ (k,n,m)=(k_j,n_j, m_j) $, where $ 3\le m_j<n_j $, for $ j=1,2 $ and $ 1\le k_1<k_2 $. Further, $ n_j\ge 4 $ for $ j=1,2 $. We return to \eqref{mukazi2} and write
\begin{eqnarray*}
\left|\Gamma_2^{(j)^{\prime}}\right|:=\left|k_j\log\rho - \log a -n_j\log\alpha\right|<\dfrac{6}{\alpha^{n_j-m_j}}, \quad \text{ for } \quad j=1,2.
\end{eqnarray*}
We do a suitable cross product between $ \Gamma_2^{(1)^{\prime}}, ~ \Gamma_2^{(2)^{\prime}} $ and $ k_1, k_2 $ to eliminate the term involving $ \log\rho $ in the above linear forms in logarithms:
\begin{eqnarray}
|\Gamma_3^{\prime}|&:=&|(k_1-k_2)\log a +(k_1n_2-k_2n_1)\log\alpha|=|k_2\Gamma_2^{(1)^{\prime}}-k_1\Gamma_2^{(2)^{\prime}}|\nonumber\\
&\le& k_2|\Gamma_2^{(1)^{\prime}}|+k_1|\Gamma_2^{(2)^{\prime}}|\quad \le \quad  \dfrac{6k_2}{\alpha^{n_1-m_1}}+\dfrac{6k_1}{\alpha^{n_2-m_2}}\quad\le \quad \dfrac{12n_2}{\alpha^{\lambda^{\prime}}},\label{basaja12}
\end{eqnarray}
where \[ \lambda^{\prime}:=\min_{1\le j\leq 2} \{n_j-m_j\}\].

We need to find an upper bound for $ \lambda^{\prime} $. If $ 12n_2/\alpha^{\lambda^{\prime}} > 1/2 $, we then get
\begin{eqnarray}
\lambda^{\prime} < \dfrac{\log (24n_2)}{\log \alpha}<4\log(24n_2).
\end{eqnarray}
Otherwise, $ |\Gamma_3^{\prime}|<\frac{1}{2} $, so
\begin{eqnarray}\label{kabaya222}
\left|e^{\Gamma_3^{\prime}}-1\right|=\left|a^{k_1-k_2}\alpha^{k_1n_2-k_2n_1}-1\right|<2|\Gamma_3^{\prime}|<\dfrac{24n_2}{\alpha^{\lambda^{\prime}}}.
\end{eqnarray}
We apply Theorem \ref{Matveev11} with the data: $ t:=2 $, $ \eta_1 := a$, $ \eta_2:= \alpha$, $ b_1:=k_1-k_2 $, $ b_2:=k_1n_2-k_2n_1 $. We take the number field $ \mathbb{K}:=\mathbb{Q}(\alpha) $ and $ D=3 $. We begin by checking that $ e^{\Gamma_3^{\prime}}-1\neq 0 $ (so $ \Gamma_3^{\prime}\neq 0 $). This is true because $ \alpha $ and $ a $ are multiplicatively independent, since $ \alpha $ is a unit in the ring of integers $ \mathbb{Q}(\alpha) $ while the norm of $ a $ is $ 1/23 $.

We note that $ |k_1-k_2|<k_2<n_2 $. Further, from \eqref{basaja12}, we have
\begin{eqnarray*}
|k_2n_1-k_1n_2|<(k_2-k_1)\dfrac{|\log a|}{\log\alpha}+\dfrac{12k_2}{\alpha^{\lambda}\log\alpha}<13k_2<13n_2
\end{eqnarray*}
given that $ \lambda\ge 1 $. So, we can take $ B:=13n_2 $.
By Theorem \ref{Matveev11}, with the same $ A_1:=\log 23 $ and $ A_2:=\log\alpha $, we have that 
\begin{eqnarray*}
\log|e^{\Gamma_3^{\prime}}-1|>-4.63\times 10^{10}(\log n_2)(\log\alpha).
\end{eqnarray*}
By comparing this with \eqref{kabaya222}, we get
\begin{eqnarray}\label{kabaya1234}
\lambda^{\prime} <1.62\times 10^{11}\log n_2.
\end{eqnarray}
Note that \eqref{kabaya1234} is better than \eqref{kabaya222}, so \eqref{kabaya1234} always holds. Without loss of generality, we can assume that $ \lambda^{\prime} = n_j-m_j $, for $ j=1,2 $ fixed.

We set $ \{j,i\}=\{1,2\} $ and return to \eqref{kabuzi22} to replace $ (k,n,m)=(k_i,n_i, m_i) $:
\begin{eqnarray}\label{muka112}
|\Gamma_1^{(i)^{\prime}}|=\left|k_i\log\rho-\log a -n_i\log\alpha-\log(1+\alpha^{m_i-n_i})\right|<\dfrac{10}{\alpha^{n_i}},
\end{eqnarray}
and also return to \eqref{mukazi2}, with $ (k, n,m)=(k_j, n_j, m_j) $:
\begin{eqnarray}\label{muka122}
|\Gamma_2^{(j)^{\prime}}|=\left|k_j\log\rho-\log a-n_j\log\alpha\right|<\dfrac{6}{\alpha^{n_j-m_j}}.
\end{eqnarray}
We perform a cross product on \eqref{muka112} and \eqref{muka122} in order to eliminate the term on $ \log\rho $:
\begin{eqnarray}
|\Gamma_4^{\prime}|&:=&\left|(k_j-k_i)\log a+(k_jn_i-k_in_j)\log\alpha+k_j\log(1+\alpha^{m_i-n_i})\right|\nonumber\\
&=&\left|k_i\Gamma_2^{(j)^{\prime}}-k_j\Gamma_1^{(i)^{\prime}}\right|\le k_i\left|\Gamma_2^{(j)^{\prime}}\right|+k_j\left|\Gamma_1^{(i)^{\prime}}\right|\nonumber\\
&<&\dfrac{6k_i}{\alpha^{n_j-m_j}}+\dfrac{10k_j}{\alpha^{n_i}}<\dfrac{16n_2}{\alpha^{\nu^{\prime}}}\label{kipro2}
\end{eqnarray}
with $ \nu^{\prime}:=\min\{n_i, n_j-m_j\} $. As before, we need to find an upper bound on $ \nu^{\prime} $. If $ 16n_2/\alpha^{\nu^{\prime}}>1/2 $, then we get
\begin{eqnarray}\label{boss122}
\nu^{\prime} < \dfrac{\log (32n_2)}{\log\alpha}< 4\log (32n_2).
\end{eqnarray}
Otherwise, $ |\Gamma_4^{\prime}|<1/2 $, so we have
\begin{eqnarray}\label{bosco122}
\left|e^{\Gamma_4^{\prime}}-1\right|\le 2|\Gamma_4^{\prime}|<\dfrac{32n_2}{\alpha^{\nu^{\prime}}}.
\end{eqnarray}
In order to apply Theorem \ref{Matveev11}, first if $ e^{\Gamma_4^{\prime}}=1 $, we obtain
\begin{eqnarray*}
a^{k_i-k_j}=\alpha^{k_jn_i-k_in_j}(1+\alpha^{-\lambda^{\prime}})^{k_j}.
\end{eqnarray*}
Since $ \alpha $ is a unit, the right-hand side in above is an algebraic integer. This is a contradiction because $ k_1<k_2 $ so $ k_i-k_j\neq 0 $, and neither $ a $ nor $ a^{-1} $ are algebraic intgers. Hence $ e^{\Gamma_4^{\prime}}\neq 1 $. By assuming that $ \nu^{\prime} \ge 100 $, we apply Theorem \ref{Matveev11} with the data:
\begin{eqnarray*}
&t:=3, \quad \eta_1:=a, \quad \eta_2:=\alpha, \quad \eta_3:=1+\alpha^{-\lambda^{\prime}},\\
& b_1:=k_j-k_i, \quad b_2:=k_jn_i-k_in_j, \quad b_3:=k_j,
\end{eqnarray*} 
and the inequalities \eqref{kabaya1234} and \eqref{bosco122}. We get
\begin{eqnarray}
\nu^{\prime}=\min\{n_i, n_j-m_j\}<1.85\times 10^{13}\lambda^{\prime}\log n_2<3\times 10^{24}(\log n_2)^{2}.
\end{eqnarray}
The above inequality also holds when $ \nu^{\prime} <100 $. Further, it also holds when the inequality \eqref{boss122} holds. So the above inequality holds in all cases. Note that the case $ \{i, j\} =\{2, 1\} $ leads to $ n_1-m_1\le n_1\le n_2+4 $ whereas $ \{i, j\} = \{1,2\} $ lead to $ \nu^{\prime} = \min\{n_1, n_2-m_2\} $. Hence, either the minimum is $  n_1 $, so
\begin{eqnarray}\label{case12}
n_1< 3\times 10^{24}(\log n_2)^{2},
\end{eqnarray}  
or the minimum is $ n_j-m_j $ and from the inequality \eqref{kabaya123} we get that
\begin{eqnarray}\label{case22}
\max_{1\le j\le 2}\{n_j- m_j\}< 3\times 10^{24}(\log n_2)^{2}.
\end{eqnarray}
Next, we assume that we are in the case \eqref{case22}. We evaluate \eqref{muka112} in $ i=1,2 $ and make a suitable cross product to eliminate the term involving $ \log\rho $:
\begin{eqnarray}
\left|\Gamma_5^{\prime}\right|&:=&\left|(k_2-k_1)\log a+(k_2n_1-k_1n_2)\log\alpha \right.\nonumber\\&&\left.+k_2\log (1+\alpha^{m_1-n_1})-k_1\log (1+\alpha^{m_2-n_2})\right|\nonumber\\
&=&\left|k_1\Gamma_1^{(2)}-k_2\Gamma_1^{(1)}\right|\le k_1\left|\Gamma_1^{(2)}\right|+k_2\left|\Gamma_1^{(1)}\right|<\dfrac{20n_2}{\alpha^{n_1}}.\label{kakawu12}
\end{eqnarray}
In the above inequality we used the inequality \eqref{faji1}to conclude that $ \min\{n_1, n_2\}\ge n_1-4 $ as well as the fact that $ n_i\ge 4 $ for $ i=1.2 $. Next, we apply a linear form in four logarithms to obtain an upper bound to $ n_1 $. As in the previous calculations, we pass from \eqref{kakawu12} to
\begin{eqnarray}\label{kakawu22}
\left|e^{\Gamma_5^{\prime}}-1\right|<\dfrac{40n_2}{\alpha^{n_1}},
\end{eqnarray}
which is implied by \eqref{kakawu12} except if $ n_1 $ is very small, say
\begin{eqnarray}\label{kakawu32}
n_1\le 4\log(40n_2).
\end{eqnarray}
Thus, we assume that \eqref{kakawu32} does not hold, therefore \eqref{kakawu22}. Then to apply Theorem \ref{Matveev11}, we fist justify that $ e^{\Gamma_5^{\prime}}\neq 1 $. Otherwise, 
\begin{eqnarray*}
a^{k_1-k_2}=\alpha^{k_2n_1-k_1n_2}(1+\alpha^{n_1-m_1})^{k_2}(1+\alpha^{n_2-m_2})^{-k_1}.
\end{eqnarray*}
By a similar argument as before, we get a contradiction. Thus, $ e^{\Gamma_{5}^{\prime}}\neq 1 $.

Then, we apply Theorem \ref{Matveev11} on the left-hand side of the inequalities \eqref{kakawu2} with the data
\begin{eqnarray*}
&t:=4, \quad \eta_1:=a, \quad \eta_2:=\alpha,\quad \eta_3:=1+\alpha^{m_1-n_1}, \quad \eta_4:=1+\alpha^{m_2-n_2},\\
&b_1:=k_2-k_1, \quad b_2:=k_2n_1-k_1n_2, \quad b_3:=k_2, \quad b_4:=k_1.
\end{eqnarray*}
Together with combining the right-hand side of \eqref{kakawu22} with the inequalities \eqref{kabaya1234} and \eqref{case22}, Theorem \ref{Matveev11} gives
\begin{eqnarray}
n_1&<&4.99\times 10^{15}(n_1-m_1)(n_2-m_2)(\log n_2)\nonumber\\
&<&2.43\times 10^{51}(\log n_2)^{4}.\label{kalo2}
\end{eqnarray}
In the above we used the facts that
\begin{eqnarray*}
\min_{1\le i\le 2}\{n_i-m_i\}<1.62\times 10^{11}\log n_2 \quad \text{and}\quad \max_{1\le i\le 2}\{n_i-m_i\}<3\times 10^{24}(\log n_2)^{2}.
\end{eqnarray*}
This was obtained under the assumption that the inequality \eqref{kakawu32} does not hold. If \eqref{kakawu32} holds, then so does \eqref{kalo2}. Thus, we have that inequality \eqref{kalo2} holds provided that inequality \eqref{case22} holds. Otherwise, inequality \eqref{case12} holds which is a better bound than \eqref{kalo2}. Hence, conclude that \eqref{kalo2} holds in all posibble cases.

By the inequality \eqref{kabizi2},
\begin{eqnarray*}
\log\rho \le k_1\log\rho \le n_1\log\alpha +\log 5 <6.92\times 10^{50}(\log n_2)^{4}.
\end{eqnarray*}
By substituting this into \eqref{lemmata1234} we get $n_2<3.67\times 10^{134}(\log n_2)^{10}$, and then, by Lemma \ref{gl}, with the data $ r:=10, ~P:=3.67\times 10^{134},~ L:=n_2 $, we get that  $ n_2< 3.07\times 10^{162}$. This immediately gives that $ n_1<4.76\times 10^{61} $.

We record what we have proved.
\begin{lemma}\label{lemmata22}
Let $(k_i, n_i, m_i)$ be a solution to $ X_{k_i}=P_{n_i}+P_{m_i} $, with $ 3\le m_i<n_i $ for $ i\in\{1,2\} $ and $ 1\le k_1<k_2 $, then
\begin{eqnarray*}
\max\{k_1, m_1\}<n_1<4.76\times 10^{61}, \quad \text{and} \quad \max\{k_2, m_2\}<n_2<3.07\times 10^{162}.
\end{eqnarray*}
\end{lemma}
\section{Reducing the bounds for $ n_1 $ and $ n_2 $ (II)}
In this section we reduce the bounds for $ n_1 $ and $ n_2 $ given in Lemma \ref{lemmata2} to cases that can be computationally treated. For this, we return to the inequalities for $ \Gamma_3^{\prime} $, $ \Gamma_4^{\prime} $ and $ \Gamma_5^{\prime} $. 
\subsection{The first reduction (II)}
We divide through both sides of the inequality \eqref{basaja12} by $(k_2-k_1)\log\alpha$. We get that
\begin{eqnarray}\label{kapa12}
\left|\dfrac{|\log a|}{\log\alpha}-\dfrac{k_2n_1-k_1n_2}{k_2-k_1}\right|<\dfrac{42n_2}{\alpha^{\lambda^{\prime}}(k_2-k_1)} \quad \text{with} \quad \lambda^{\prime}: = \min_{1 \le i \le 2}\{n_i-m_i\}.
\end{eqnarray}
We assume that $ \lambda^{\prime}\ge 10 $. Below we apply Lemma \ref{Legendre}. We put $ \tau^{\prime}:= \frac{|\log a|}{\log\alpha} $, which is irrational and compute its continued fraction \[ [a_0, a_1, a_2, \ldots]=[1, 6, 2, 1, 18, 166, 1, 2, 13, 1, 2, 5, 1, 5, 1, 2, 3, 1, 1, 31, 1, 3, 2, 3, \ldots] \] and its convergents \[ \left[\frac{p_0}{q_0},\frac{p_1}{q_1}, \frac{p_2}{q_2}, \ldots\right] =\left[ 1, \frac{7}{6}, \frac{15}{13}, \frac{22}{19}, \frac{411}{355}, \frac{68248}{58949}, \frac{68659}{59304}, \
\frac{205566}{177557}, \frac{2741017}{2367545}, \frac{2946583}{2545102},\ldots \right]. \]
Furthermore, we note that taking $ N:=3.07\times 10^{162} $ (by Lemma \ref{lemmata22}), it follows that
\begin{eqnarray*}
q_{296}>N>n_2>k_2-k_1 \quad \text{and}\quad a(N):=\max\{a_j:0\le j\le 296\}=a_{189}=1028.
\end{eqnarray*}
Thus, by Lemma \ref{Legendre}, we have that
\begin{eqnarray}\label{kapa22}
\left|\tau^{\prime} - \dfrac{k_2n_1-k_1n_2}{k_2-k_1}\right|>\dfrac{1}{1030(k_2-k_1)^{2}}.
\end{eqnarray}
Hence, combining the inequalities \eqref{kapa12} and \eqref{kapa22}, we obtain
\begin{eqnarray*}
\alpha^{\lambda^{\prime}}<43260n_2(k_2-k_1)<4.08\times 10^{329},
\end{eqnarray*}
so $ \lambda^{\prime}\le 2661 $. This was obtained under the assumption that $ \lambda^{\prime}\ge 10 $, Otherwise, $ \lambda^{\prime}<10<2661 $ holds as well.

Now, for each $ n_i-m_i = \lambda^{\prime}\in [1, 2661] $ we estimate a lower bound $ |\Gamma_4^{\prime}| $, with
\begin{eqnarray}\label{LLL12}
\Gamma_4^{\prime}&=&(k_j-k_i)\log a+(k_jn_i-k_in_j)\log\alpha+k_j\log(1+\alpha^{m_i-n_i})
\end{eqnarray}
given in the inequality \ref{kipro2}, via the same procedure described in Subsection \ref{Reduction}  (LLL-algorithm). We recall that $ \Gamma_4^{\prime}\neq 0 $. 

We apply Lemma \ref{LLL} with the data: 
\begin{eqnarray*}
&t:=3, \quad \tau_1:=\log a, \quad \tau_2:=\log\alpha, \quad \tau_3:=\log(1+\alpha^{-\lambda^{\prime}}),\\
&x_1:=k_j-k_i, \quad x_2:=k_jn_i-k_in_j, \quad x_3:=k_j.
\end{eqnarray*}
We set $ X:= 3.99\times 10^{163} $ as an upper bound to $ |x_i|<13n_2 $ for all $ i=1, 2, 3 $, and $ C:=(20X)^{5} $. A computer in \textit{Mathematica} search allows us to conclude, together with the inequality \eqref{kipro2}, that
\begin{eqnarray*}
8\times 10^{-660}<\min_{1\le \lambda \le  2661}|\Gamma_4^{\prime}|<16n_2\alpha^{-\nu^{\prime}},\quad \text{with} \quad \nu^{\prime}:=\min\{n_i, n_j-m_j\}
\end{eqnarray*}
which leads to $ \nu^{\prime} \le 6643 $. As we have noted before, $ \nu^{\prime} = n_1 $ (so $n_1\le 6643$) or $ \nu^{\prime} = n_j-m_j $.

Next, we suppose that $ n_j-m_j = \nu^{\prime} \le 6643 $. Since $ \lambda^{\prime}\le 2661 $, we have
\begin{eqnarray*}
\lambda^{\prime} := \min_{1\le i \le 2}\{n_i-m_i\}\le 2661 \quad \text{and} \quad \chi^{\prime} :=\max_{1\le i \le 2}\{n_i-m_i\}\le 6643.
\end{eqnarray*}
Now, returning to the inequality \eqref{kakawu12} which involves
\begin{eqnarray}\label{LLL122}
\Gamma_5^{\prime}:&=&(k_2-k_1)\log a+(k_2n_1-k_1n_2)\log\alpha \nonumber\\&&+k_2\log (1+\alpha^{m_1-n_1})-k_1\log (1+\alpha^{m_2-n_2})\neq 0,
\end{eqnarray}
we use again the LLL-algorithm to estimate the lower bound for $ |\Gamma_5^{\prime}| $ and thus, find a bound for $ n_1 $ that is better than the one given in Lemma \ref{lemmata22}.

We distinguish the cases $ \lambda^{\prime} < \chi^{\prime} $ and $ \lambda^{\prime} = \chi^{\prime} $. 
\subsection{The case $ \lambda^{\prime} < \chi^{\prime} $}
We take $ \lambda^{\prime} \in [1, 2661] $ and $ \chi^{\prime} \in [\lambda^{\prime}+1, 6643] $ and apply Lemma \ref{LLL} with the data:
\begin{eqnarray*}
& t:=4,\quad \tau_1:=\log a, \quad \tau_2:= \log\alpha, \quad \tau_3: = \log(1+\alpha^{m_1-n_1}), \quad  \tau_4: = \log(1+\alpha^{m_2-n_2}),\\
& x_1:=k_2-k_1, \quad x_2:= k_2n_1-k_1n_2, \quad x_3: = k_2, \quad x_4:=-k_1.
\end{eqnarray*}
We also put $ X:=3.99\times 10^{163} $ and $ C:=(20X)^{9} $. As before, after a computer search in \textit{Mathematica} together with the inequality \ref{kakawu12}, we can confirm that
\begin{eqnarray}
9.9\times 10^{-1317}<\min_{\substack{1\le \lambda\le 2661 \\ \lambda+1\le \chi \le 6643}}|\Gamma_5^{\prime}| < 20n_2 \alpha^{-n_1}.
\end{eqnarray}
This leads to the inequality
\begin{eqnarray}
\alpha^{n_1}< 2.02\times 10^{1317}n_2.
\end{eqnarray}
Subsitituting for the bound  $ n_2 $ given in Lemma \ref{lemmata22}, we get that $ n_1\le 11948 $.
\subsection{The case$ \lambda^{\prime} =\chi^{\prime} $}
In this case, we have 
\begin{eqnarray*}
\Lambda_5^{\prime}:=(k_2-k_1)(\log a+\log(1+\alpha^{m_1-n_1}))+(k_2n_1-k_1n_2)\log\alpha \neq 0.
\end{eqnarray*}
We divide through the inequality \ref{kakawu12} by $(k_2-k_1)\log\alpha$ to obtain
\begin{eqnarray}\label{kaka2}
\left|\dfrac{|\log a+\log(1+\alpha^{m_1-n_1})|}{\log \alpha}-\dfrac{k_2n_1-k_1n_2}{k_2-k_1}\right|<\dfrac{70n_2}{\alpha^{n_1}(k_2-k_1)}
\end{eqnarray}
We now put $$ \tau_{\lambda^{\prime}}:=\dfrac{|\log a+\log(1+\alpha^{-\lambda^{\prime}})|}{\log \alpha} $$ and compute its continued fractions $ [a_0^{(\lambda^{\prime})}, a_1^{(\lambda^{\prime})} , \ldots] $ and its convergents $[p_0^{(\lambda)}/q_0^{(\lambda^{\prime})}, p_1^{(\lambda^{\prime})}/q_1^{(\lambda^{\prime})}, \ldots]$ for each $ \lambda^{\prime} \in [1, 2661] $. Furthermore, for each case we find an integer $ t_{\lambda^{\prime}} $ such that $ q_{t_{\lambda^{\prime}}}^{(\lambda^{\prime})}>N:=3.07\times 10^{162}>n_2>k_2-k_1 $ and calculate
\begin{eqnarray*}
a(N):=\max_{1\le\lambda^{\prime}\le 2661}\left\{a_{i}^{(\lambda^{\prime})}: 0 \le i \le t_{\lambda^{\prime}}\right\}.
\end{eqnarray*}
A computer search in \textit{Mathematica} reveals that for  $ \lambda^{\prime} = 2466  $, $ t_{\lambda^{\prime}} = 298 $ and $ i= 295 $, we have that $ a(N) = a_{295}^{(2466)}=2818130 $. Hence, combining the conclusion of Lemma \ref{Legendre} and the inequality \eqref{kaka2}, we get
\begin{eqnarray*}
\alpha^{n_1}< 70\times 2818132n_2(k_2-k_1)< 1.86\times 10^{333},
\end{eqnarray*}
so $ n_1\le 2690 $. Hence, we obtain that $ n_1\le 11948 $ holds in all cases ($\nu^{\prime} =n_1$, $\lambda^{\prime} < \chi^{\prime}$ or $\lambda^{\prime} = \chi^{\prime} $).  By the inequality \eqref{kabizi2}, we have that
\begin{eqnarray*}
\log\rho\le k_1\log\rho \le n_1\log\alpha +\log 5 <3410.
\end{eqnarray*}
By considering the second inequality in \eqref{lemmata1234}, we can conclude that $ n_2\le 8.17\times 10^{39}(\log n_2)^2 $, which yields $ n_2<2.76\times 10^{44}  $, by a simple application of Lemma \ref{gl} as before. Below, we summarise the first cycle of our reduction process:
\begin{eqnarray}
n_1\le 11948 \quad \text{and} \quad n_2\le 2.76\times 10^{44}.
\end{eqnarray}
As in the previous case, from the above, we note that the upper bound on $ n_2 $ represents a very good reduction of the bound given in Lemma \ref{lemmata22}. Hence, we expect that if we restart our reduction cycle with the new bound on $ n_2 $, then we get a better bound on $ n_1 $. Thus, we return to the inequality \eqref{kapa2} and take $ N:=2.76\times 10^{44} $. A computer search in \textit{Mathematica} reveals that
\begin{eqnarray*}
q_{88}>N>n_2>k_2-k_1 \quad \text{and} \quad a(N):=\max\{a_i: 0\le i\le 88\}=a_{55} =397,
\end{eqnarray*}
from which it follows that $ \lambda \le 738 $. We now return to \eqref{LLL12} and we put $ X:=2.76\times 10^{44} $ and $ C:=(10X)^{5} $ and then apply the LLL-algorithm in Lemma \ref{LLL} to $ \lambda \in[1, 738] $. After a computer search, we get
\begin{eqnarray*}
8.6\times 10^{-183}<\min_{1\le \lambda^{\prime}\le 738}|\Gamma_4^{\prime}| < 16n_2\alpha^{-\nu^{\prime}},
\end{eqnarray*}
then $ \nu^{\prime} \le 1838 $. By continuing under the assumption that $ n_j-m_j=\nu \le 1838 $, we return to \eqref{LLL122} and put $ X:=2.76\times 10^{44} $, $ C:=(10X)^9 $ and $ N:=2.76\times 10^{44} $ for the case $ \lambda^{\prime} <\chi^{\prime} $ and $ \lambda^{\prime} = \chi^{\prime} $. After a computer search, we confirm that
\begin{eqnarray*}
8\times 10^{-365}< \min_{\substack{1\le \lambda \le 738\\ \lambda+1\le \chi \le 1838}}|\Gamma_5^{\prime}|<6n_2\alpha^{-n_1},
\end{eqnarray*}
gives $ n_1\le 3304 $, and
$
a(N)=a_{125}^{(160)}=155013
$, leads to $ n_1\le 774 $. Hence, in both cases $ n_1\le 3304 $ holds. This gives $ n_2\le 4\times 10^{42} $ by a similar procedure as before, and $ k_1\le  $.

We record what we have proved.

\begin{lemma}\label{firstredn2}
Let $ (k_i, n_i, m_i) $ be a solution to $ X_i=P_{n_i}+P_{m_i} $, with $ 3\le m_i<n_i $ for $ i=1,2 $ and $ 1\le k_1<k_2 $, then
\begin{eqnarray*}
m_1<n_1\le 3304, \quad  k_1\le 3108 \quad \text{and} \quad n_2\le 4\times 10^{42}.
\end{eqnarray*}
\end{lemma}
\subsection{The final reduction (II)}
Returning back to \eqref{kay2} and \eqref{kayiii2} and using the fact that $ (X_1, X_1) $ is the smallest positive solution to the Pell equation \eqref{Pelleqn2}, we obtain
\begin{eqnarray*}
X_k &=& \rho^k+\varrho^k \quad =\quad \left(\dfrac{(X_1+Y_1\sqrt{d}}{2}\right)^{k}+\left(\dfrac{X_1-Y_1\sqrt{d}}{2}\right)^{k}\\
&=&\left(\dfrac{X_1+\sqrt{X_1^2\mp 4}}{2}\right)^{k}+\left(\dfrac{X_1-\sqrt{X_1^2\mp 4}}{2}\right)^{k}: = R^{\pm}_{k}(X_1).
\end{eqnarray*}
Thus, we return to the Diophantine equation $ X_{k_1}=P_{n_1}+P_{m_1} $ and consider the equations 
\begin{eqnarray}\label{water2}
R^{+}_{k_1}(X_1)=P_{n_1}+P_{m_1} \quad \text{and} \quad R^{-}_{k_1}(X_1)=P_{n_1}+P_{m_1},
\end{eqnarray}
with $ k_1\in [1, 3108] $, $ m_1\in [3,3304 ] $ and $n_1\in[m_1+1, 3304 ] $.

A computer search in \textit{Mathematica} on the above equations in \eqref{water2} shows that there are only finitely many solutions that we list in the tables below. We note that $$ \dfrac{3+\sqrt{5}}{2}=\left(\dfrac{1+\sqrt{5}}{2}\right)^2  \qquad \text{and}\qquad 2+\sqrt{5} = \left(\dfrac{1+\sqrt{5}}{2}\right)^3, $$
so these come from the same Pell equation with $ d=5 $. Similarly,
\begin{eqnarray*}
\dfrac{11+\sqrt{13}}{2} = \left(\dfrac{3+\sqrt{13}}{2}\right)^2, \qquad \text{and} \qquad \dfrac{51+7\sqrt{53}}{2}=\left(\dfrac{7+\sqrt{53}}{2}\right)^{2}
\end{eqnarray*}
these also come from the same Pell equation with $ d=13 $ and $ d=53 $, respectively.
\begin{table}[H]
\parbox{.38\linewidth}{
\centering
\begin{tabular}{|c c c c c|}
\hline
& & &$ R^{+}_{k_1}(X_1) $&\\
\hline
$ k_1 $& $ X_1 $& $ Y_1 $& $ d $ & $ \rho $\\
\hline
$ 2 $&$ 3 $& $ 1 $&$ 5 $ & $ (3+\sqrt{5})/2 $\\
$ 2 $&$ 4 $& $ 2 $&$ 3 $ & $ 2+\sqrt{3} $\\
$ 2 $&$ 5 $& $ 1 $&$ 21 $ & $ (5+\sqrt{21})/2 $\\
$ 3 $&$ 9 $& $ 1 $&$ 77 $ & $ (9+\sqrt{77})/2 $\\
$ 2 $&$ 10 $& $ 4 $&$ 6 $ & $ 5+ 2\sqrt{6}$\\
$ 2 $&$ 11 $& $ 3 $&$ 13 $ & $ (11+3\sqrt{13})/2 $\\
$ 2 $&$ 12 $& $ 2 $&$ 35 $ & $ 6+\sqrt{35} $\\
$ 2 $&$ 13 $& $ 1 $&$ 165 $ & $ (13+2\sqrt{165})/2 $\\
$ 3 $&$ 15 $& $ 1 $&$ 221 $ & $ (15+\sqrt{221})/2 $\\
$ 2 $&$ 25 $& $ 3 $&$ 69 $ & $ (25+3\sqrt{69})/2 $\\
$ 2 $&$ 44 $& $ 2 $&$ 483 $ & $ 22+\sqrt{483} $\\
$ 2 $&$ 51 $& $ 7 $&$ 53 $ & $ (51+7\sqrt{53})/2 $\\
$ 2 $&$ 88 $& $ 6 $&$ 215 $ & $ 44+3\sqrt{215} $\\
$ 2 $&$ 2570 $& $ 4 $&$ 412806 $ & $ 1285+2\sqrt{412806} $\\
\hline
\end{tabular}
}
\hfill
\parbox{.40\linewidth}{
\centering
\begin{tabular}{|c c c c c|}
\hline
& & &$ R^{-}_{k_1}(X_1) $&\\
\hline
$ k_1 $& $ X_1 $& $ Y_1 $& $ d $ & $ \rho $\\
\hline
$ 2 $&$ 1 $& $ 1 $&$ 5 $ & $ (1+\sqrt{5})/2 $\\
$ 2 $&$ 2 $& $ 2 $&$ 2 $ & $ 1+\sqrt{2} $\\
$ 2 $&$ 3 $& $ 1 $&$ 13 $ & $ (3+\sqrt{13})/2 $\\
$ 2 $&$ 4 $& $ 2 $&$ 5 $ & $ 2+\sqrt{5} $\\
$ 2 $&$ 6 $& $ 2 $&$ 10 $ & $ 3+\sqrt{10} $\\
$ 2 $&$ 7 $& $ 1 $&$ 53 $ & $ (7+\sqrt{53})/2 $\\
$ 2 $&$ 8 $& $ 2 $&$ 17 $ & $ 4+\sqrt{17} $\\
$ 2 $&$ 10 $& $ 2 $&$ 26 $ & $ 5+\sqrt{26} $\\
$ 2 $&$ 11 $& $ 5 $&$ 5 $ & $ (11+5\sqrt{5})/2 $\\
$ 2 $&$ 19 $& $ 1 $&$ 365 $ & $ (19+\sqrt{365})/2 $\\
$ 2 $&$ 22 $& $ 2 $&$ 122 $ & $ 11+\sqrt{122} $\\
$ 2 $&$ 30 $& $ 2 $&$ 226 $ & $ 15+\sqrt{226} $\\
$ 2 $&$ 58 $& $ 2 $&$ 842 $ & $ 29+\sqrt{842} $\\
$ 2 $&$ 88 $& $ 2 $&$ 1937 $ & $ 44+\sqrt{1937} $\\
$ 2 $&$ 178 $& $ 2 $&$ 7922 $ & $ 89+\sqrt{7922} $\\
$ 2 $&$ 3480 $& $ 2 $&$3027601 $ & $ 1740+\sqrt{3027601} $\\
\hline
\end{tabular}
}
\end{table}

From the above tables, we set each $ \rho:=\rho_{t} $ for $ t=1, 2, \ldots 25 $.
We then work on the linear forms in logarithms $ \Gamma_1^{\prime} $ and $ \Gamma_2^{\prime} $, in order to reduce the bound on $ n_2 $ given in  Lemma \ref{firstredn2}. From the inequality \eqref{mukazi2}, for $ (k,n,m):=(k_2, n_2, m_2) $, we write
\begin{eqnarray}\label{mukazi1112}
\left|k_2\dfrac{\log\rho_t}{\log\alpha}-n_2+\dfrac{\log a}{\log (\alpha^{-1})}\right|<\left(\frac{6}{\log \alpha}\right)\alpha^{-(n_2-m_2)}, \quad \text{for}\quad t=1,2, \ldots 25.
\end{eqnarray}
We put
\begin{eqnarray*}
\tau_{t}:=\dfrac{\log\rho_t}{\log\alpha}, \qquad \mu_t:=\dfrac{\log a}{\log(\alpha^{-1})}\qquad \text{and} \quad (A_t, B_t):=\left(\frac{6}{\log \alpha}, \alpha\right).
\end{eqnarray*}
We note that $ \tau_t $ is transcendental by the Gelfond-Schneider's Theorem and thus, $ \tau_t $ is irrational. We can rewrite the above inequality, \ref{mukazi1112} as
\begin{eqnarray}\label{mukaziii112}
0<|k_2\tau_t-n_2+\mu_t|<A_tB_t^{-(n_2-m_2)}, \quad \text{for} \quad t=1, 2, \ldots, 25.
\end{eqnarray}
We take $ N:= 4\times 10^{42} $ which is the upper bound on $ n_2 $ according to Lemma \ref{firstredn2} and apply Lemma \ref{Dujjella} to the inequality \eqref{mukaziii112}. As before, for each $ \tau_t $ with $ t=1, 2, \ldots, 26 $, we compute its continued fraction $ [a_0^{(t)}, a_1^{(t)}, a_2^{(t)}, \ldots ] $ and its convergents $ p_0^{(t)}/q_0^{(t)}, p_1^{(t)}/q_1^{(t)}, p_2^{(t)}/q_2^{(t)}, \ldots $. For each case, by means of a computer search in \textit{Mathematica}, we find and integer $s_{t}$ such that
\begin{eqnarray*}
q^{(t)}_{s_t}> 2.4\times 10^{43}=6N \qquad \text{ and } \qquad \epsilon_t:=||\mu_t q^{(t)}||-N||\tau_t q^{(t)}|>0.
\end{eqnarray*}
We finally compute all the values of $ b_t:=\lfloor \log(A_t q^{(t)}_{s_t}/\epsilon_t)/\log B_t \rfloor $. The values of  $ b_t $ correspond to the upper bounds on $ n_2-m_2 $, for each $ t=1, 2, \ldots, 25 $, according to Lemma \ref{Dujjella}. We record the results of the computations for each $ t $ in the following table.

\begin{center}
\begin{tabular}{lllccc}
\hline 
$t$ & $\rho_t$ & $s_t$ & $q_{s_t} $ & $\epsilon_t>$ & $b_t$ \\ 
\hline 
$1$ & $1+\sqrt{2}$ & $78$ & $ 1.46195\times 10^{44} $ & $ 0.1578 $ & $379$ \\ 
$2$ & $ 2+\sqrt{3} $ & 100 & $ 8.93366\times 10^{43} $ & $ 0.3147 $ & $374$ \\ 
$3$ & $ (1+\sqrt{5})/2 $ & $ 82 $ & $ 2.96985\times 10^{43} $ & $ 0.4479 $ & $369$ \\  
$4$ & $ 5+2\sqrt{6} $ & $ 80 $ & $ 3.16032\times 10^{43} $ & $ 0.1940 $ & $372$ \\  
$5$ & $ 3+\sqrt{10} $& $ 88 $ & $ 4.21425\times 10^{43} $ & $ 0.2358 $ & $373$ \\ 
$6$ & $(3+\sqrt{13})/2$ & $ 91 $ & $6.62314\times 10^{43}$ & $0.0666$ & $379$ \\  
$7$ & $4+\sqrt{17} $ & $ 92 $ & $ 1.11753\times 10^{44} $ & $ 0.2387 $ & $376$ \\  
$8$ & $(5+\sqrt{21})/2 $ & $ 73 $ & $  2.44965 \times 10^{43}$ & $ 0.0400 $ & $377$ \\  
$9$ & $5+\sqrt{26} $ & $ 98 $ & $ 3.23107\times 10^{43} $ & $ 0.2333 $ & $372$ \\ 
$10$ & $6+\sqrt{35}$  & $ 83 $ & $ 1.87425\times 10^{44} $ & $ 0.1172 $& $381$ \\  
$11$ & $(7+\sqrt{53})/2$ & $ 96 $ & $ 1.82440\times 10^{44} $ & $ 0.3875 $ & $376$ \\ 
$12$ & $(25+3\sqrt{69})/2$ & $ 80 $ & $ 2.40911\times 10^{43} $ & $ 0.2013 $ & $371$ \\  
$13$ & $ (9+\sqrt{77})/2 $ & $ 82 $ & $ 2.54747\times 10^{43} $ & $ 0.1470 $ & $373$ \\  
$14$ & $11+\sqrt{122} $ & $ 76 $ & $ 4.91937\times 10^{44} $ & $ 0.4004 $ & $380$ \\  
$15$ & $(13+2\sqrt{165})/2 $ & $ 86 $ & $ 2.61323\times 10^{43} $ & $ 0.1664 $ & $372$ \\  
$16$ & $44+3\sqrt{215} $ & $ 80 $ & $ 3.14146\times 10^{43} $ & $ 0.3298 $ & $371$ \\  
$17$ & $(15+\sqrt{221})/2$ & $75$ & $5.70467\times 10^{43}$ & $ 0.4661 $ & $371$ \\ 
$18$ & $15+\sqrt{226}$  & $ 79 $ & $ 4.78438\times 10^{43} $ & $ 0.4046 $& $371$ \\  
$19$ & $(19+\sqrt{365})/2$ & $ 78 $ & $ 3.05270\times 10^{43} $ & $ 0.1985 $ & $372$ \\ 
$20$ & $ 22+\sqrt{483} $ & $ 75 $ & $ 4.19689\times 10^{43} $ & $ 0.1559 $ & $374$ \\  
$21$ & $29+\sqrt{842} $ & $ 87 $ & $ 8.14707\times 10^{44} $ & $ 0.2964 $ & $382$ \\  
$22$ & $44+\sqrt{1937}$ & $ 87 $ & $ 4.70884\times 10^{43} $ & $ 0.1191 $ & $376$ \\  
$23$ & $89+\sqrt{7922} $ & $ 79 $ & $ 2.43413\times 10^{43} $ & $ 0.4418 $ & $369$ \\  
$24$ & $1285+2\sqrt{412806} $ & $85 $ & $ 2.22078\times 10^{45} $ & $ 0.4501 $ & $385$ \\ 
$25$ & $1740+\sqrt{3027601}$ & $77$ & $2.33761\times 10^{44}$ & $ 0.3352 $ & $378$ \\ 
\hline 
\end{tabular} 
\end{center}
\vspace*{0.5cm}

By replacing $ (k, n, m):=(k_2, n_2, m_2) $ in the inequality \eqref{kabuzi22}, we can write
\begin{eqnarray}\label{water1122}
\left|k_2\dfrac{\log\delta_t}{\log\alpha}-n_2+\dfrac{\log(a(1+\alpha^{-(n_2-m_2)}))}{\log(\alpha^{-1})}\right|< \left(\dfrac{10}{\log\alpha}\right)\alpha^{-n_2}, 
\end{eqnarray}
for  $t=1,2,\ldots, 25$.
We now put
\begin{eqnarray*}
\tau_{t}:=\dfrac{\log\delta_t}{\log\alpha}, \quad \mu_{t, n_2-m_2}:=\dfrac{\log(a(1+\alpha^{-(n_2-m_2)}))}{\log(\alpha^{-1})}\quad \text{and} \quad (A_t, B_t):=\left(\frac{10}{\log \alpha}, \alpha\right).
\end{eqnarray*}
With the above notations, we can rewrite \eqref{water1122} as
\begin{eqnarray}\label{water212}
0<|k_2\tau_t - n_2+\mu_{t, n_2-m_2}|<A_tB_t^{-n_2}, \quad \text{ for} \quad t=1,2, \ldots 25.
\end{eqnarray}
We again apply Lemma \ref{Dujjella} to the above inequality \eqref{water212}, for 
\begin{eqnarray*}
t=1, 2, \ldots, 25, \quad n_2-m_2 =1, 2, \ldots, b_t, \quad \text{with}\quad N:=4\times 10^{43}.
\end{eqnarray*}
We take
\begin{eqnarray*}
\epsilon_{t, n_2-m_2}:=||\mu_t q^{(t, n_2-m_2)}|| -N||\tau_t q^{(t, n_2-m_2)}||>0,
\end{eqnarray*}
and 
\begin{eqnarray*}
b_{t, n_2-m_2}:=\lfloor \log(A_t q^{(t, n_2-m_2)}_{s_t}/\epsilon_{t, n_2-m_2})/\log B_t \rfloor.
\end{eqnarray*}
With the help of \text{Mathematica}, we obtain that

\vspace{0.4cm}
\begin{tabular}{c|cccccccccccccc}
\hline
$t$& $ 1 $&$ 2 $&$ 3 $&$ 4 $&$ 5 $&$ 6 $&$ 7 $&$ 8 $&$ 9 $&$ 10 $&$ 11 $&$ 12 $&$ 13 $\\

$b_{t, n_2-m_2}$& $ 398 $&$ 404 $&$ 399 $&$ 413 $&$ 390 $&$ 398 $&$ 401 $&$ 397 $&$ 390$&$ 413 $&$ 401 $&$ 396 $&$ 396 $&\\
\hline
$t$&$ 14 $&$ 15 $&$ 16 $&$ 17 $&$ 18 $&$ 19 $ &$ 20 $&$ 21 $&$ 22 $&$ 23 $&$ 24 $&$ 25 $\\

$b_{t, n_2-m_2}$&$ 402 $& $393 $&$ 395 $&$ 392 $&$ 401 $&$ 396 $&$ 392 $& $ 400 $&$ 401 $&$ 392 $&$ 414 $&$ 395$\\
\hline
\end{tabular}

\begin{eqnarray*}
\max\{b_{t, n_2-m_2}: t=1, 2, \ldots, 25 \quad \text{and} \quad n_2-m_2 = 1, 2, \ldots d_t\} \le 414.
\end{eqnarray*}
Thus, by Lemma \ref{Dujjella}, we have that $ n_2\le 414 $, for all $ t=1,2, \ldots, 25 $, and by the inequality \eqref{faji2} we also have that $ n_1\le n_2+4 $. From the fact that $ \rho^{k}\le 2\alpha^{n+1} $, we can conclude that $ k_1< k_2\le 248  $. Collecting everything together, our problem is reduced to search for the solutions for \eqref{kalayim12} in the following range
\begin{eqnarray}
\quad 1\le k_1<k_2\le 248, \quad 0\le m_1<n_1 \in [3, 414] \quad \text{and} \quad 0\le m_2<n_2 \in [3, 414].
\end{eqnarray}
After a computer search on the equation \eqref{kalayim12} on the above ranges, we obtained the following solutions, which are the only solutions for the  exceptional $ d $ cases  we have stated in Theorem \ref{Main2}:

For the  $ +4 $ case:
\begin{eqnarray*}
(d=3)&&X_1=4=P_7+P_0=P_6+P_3=P_5+P_5, \quad X_2=14=P_{11}+P_5=P_{10}+P_8,\\ && X_3=52=P_{16}+P_6;\\
(d=5)&&X_1=3=P_6+P_0 = P_5+P_3, \quad  X_2=7=P_9+P_0=P_7+P_6,~ \\&&  X_3=18=P_{12}+P_5;\\
  (d=21)&&X_1 =5=P_8+P_0 = P_7+P_3= P_6+P_5,\quad X_2=23=P_{13}+P_5=P_{12}+P_{9}, \\
&&  X_3=2525=P_{30}+P_{11}.
\end{eqnarray*}

For the $ -4 $ case:
\begin{eqnarray*}
(d=2) && X_1=2=P_5+P_0=P_3+P_3, \quad X_2=14=P_{11}+P_5=P_{10}+P_{8};\\
(d=5)&& X_1=1=P_3+P_0, \quad X_2=4=P_7+P_0=P_6+P_3=P_5+P_5, \\ && X_3=11=P_{10}+P_5=P_{9}+P_7, \quad X_4=29=P_{14}+P_{3}.
\end{eqnarray*}
This completes the proof of Theorem \ref{Main2}. \qed

\section*{Acknowledgements}
The author was supported by the Austrian Science Fund (FWF) grants: F5510-N26 -- Part of the special research program (SFB), ``Quasi Monte Carlo Metods: Theory and Applications'' and W1230 --``Doctoral Program Discrete Mathematics''. Part of this paper was written when the author visited the Institut de Math\'ematiques de Bordeaux, Universit\'e de Bordeaux, in May 2019. He would like to thank this institution for its hospitality and the fruitful working environment.

 \end{document}